\documentclass[11pt]{article}

\usepackage{exscale,amssymb}
 \usepackage{graphicx}

\usepackage{comment}
\usepackage{amssymb}
\usepackage{amsmath}
\usepackage{amsthm}
\usepackage{mathtools}
\usepackage{algorithmicx}
 \usepackage[ruled]{algorithm}
 \usepackage{algpseudocode}
 \usepackage{algpascal}
 \usepackage{algc}
 \usepackage{color}
 \usepackage[color]{changebar}
  \alglanguage{pseudocode}

\newtheorem{theorem}{Theorem}%[section]
\newtheorem{remark}{Remark}%[section]

\def\R{\mathbb{R}}
\def\Z{\mathbb{Z}}

\def\dfrac#1#2{{\displaystyle{#1\over#2}}}

\let\oldref\ref
\def\ref#1{{\normalfont\oldref{#1}}}
\def\eqref#1{{\normalfont(\oldref{#1})}}

\mathchardef\Gamma="7100
\mathchardef\Delta="7101
\mathchardef\Theta="7102
\mathchardef\Lambda="7103
\mathchardef\Xi="7104
\mathchardef\Pi="7105
\mathchardef\Sigma="7106
\mathchardef\Upsilon="7107
\mathchardef\Phi="7108
\mathchardef\Psi="7109
\mathchardef\Omega="710A

\definecolor{cmg}{RGB}{250, 0, 250}
\definecolor{cmm}{RGB}{150, 10, 210}

\title{A low-cost alternating projection approach for a continuous formulation of convex and cardinality constrained optimization}

\author{N. Kreji\'c\thanks{Department of Mathematics and Informatics, Faculty of Sciences, University of Novi Sad, Trg Dositeja Obradovi\'ca 4, 21000 Novi Sad, Serbia (natasak@uns.ac.rs).	}
 \and E. H. M. Krulikovski\thanks{Center for Mathematics and Applications (NovaMath), FCT NOVA, 2829-516
	 Caparica, Portugal (e.krulikovski@fct.unl.pt). }
\and M. Raydan\thanks{Center for Mathematics and Applications (NovaMath), FCT NOVA, 2829-516 Caparica, Portugal (m.raydan@fct.unl.pt).}
}

\cbcolor{blue}

\begin{document}

\date{September 6, 2022}

\maketitle

\begin{abstract}
We consider convex constrained optimization problems that also include a cardinality constraint.  In general, optimization problems with cardinality
constraints are  difficult mathematical programs which are usually solved by global techniques from discrete optimization.
We assume that the region defined by the convex constraints can be written as the intersection of a finite collection of convex sets,
 such that it is easy and inexpensive to project onto each one of them (e.g., boxes, hyper-planes, or half-spaces).
 Taking advantage of a recently developed continuous reformulation that relaxes the cardinality constraint, we propose a specialized penalty
 gradient projection scheme combined with alternating projection ideas to solve these problems. To illustrate the combined scheme, we focus on
 the standard mean-variance portfolio optimization problem for which we can only invest in a preestablished limited number of assets.
 For these portfolio problems with cardinality constraints we present a numerical study on a variety of data sets involving real-world capital
 market indices from major stock markets.  On those data sets we illustrate the computational performance of the proposed scheme to produce the
  effective frontiers for different values of the limited number of allowed assets.     \\ [2mm]
{\bf AMS Subject Classification:}  90C30, 65K05, 91G10, 91G15.  \\ [2mm]
{\bf Keywords:}  Cardinality constraints, Portfolio optimization, Efficient frontier, Projected gradient methods, Dykstra's algorithm.
\end{abstract}

\section{Introduction} \label{intro}

We are interested in convex constrained optimization problems with an additional cardinality constraint.
In other words, we are interested in finding sparse solutions of those optimization problems, i.e.
solutions with a limited number of nonzero elements, as required in many areas including
image and signal processing, mathematical statistics, machine learning, portfolio optimization problems,
 among others. One effective way  to ensure the sparsity of the obtained solution is imposing a cardinality constraint
  where the number of nonzero elements of the solution is bounded in advance. \\

To be precise,  let us consider the following constrained optimization problem
\begin{equation} \label{origprb}
	\displaystyle\min_{x} f(x) \;\;\;\; \mbox{ subject to }\;\;\;\;  x\in\Omega \; \mbox{ and } \; \|x\|_0 \leq \alpha,
\end{equation}
where $f:\mathbb{R}^{n}\to \mathbb{R}$ is  continuously differentiable, $1\leq \alpha < n$ is a given
natural number, $\Omega$  is a convex subset of $\mathbb{R}^{n}$  (that will change depending on the considered application),
and  the  $L_0$ (quasi) norm $\|x\|_0$ denotes the number of nonzero components of $x$.  The sparsity constraint $\|x\|_0\leq \alpha$
is also called the cardinality constraint. Of course,  we will assume that $\alpha < n$ since otherwise the cardinality constraint could be discarded. \\

The main difference between problem (\ref{origprb}) and a standard convex constrained optimization problem is that the cardinality constraint,
despite of the notation, is not a norm, nor continuous neither convex.
 Because of the non-tractability of the  so-called zero norm  $\|x\|_0$, the 1-norm $\|x\|_1$  has also been frequently
 considered to  develop good approximate algorithms.
 Clearly, to impose a required level of sparsity, the use of the zero norm in (\ref{origprb}) is much more effective. \\

 Optimization problems with cardinality constraints are very hard problems which are typically solved by global
 techniques from discrete or combinatorial optimization.  However, in a more general setting, a continuous reformulation has been recently
 proposed and analyzed in \cite{Burdakov} to deal with this difficult
  cardinality constraint. The main idea is to address the continuous counterpart of problem (\ref{origprb}):
  \begin{equation} \label{relaxed}
  	 \begin{array}{ll}
\displaystyle	\min_{x,y}   & f(x) \\
	\mbox{subject to:} &  x\in\Omega,  \\
       	& e^{\top} y \geq n-\alpha,  \\
      	& x_i y_i = 0, \; \mbox{ for all } 1\leq i \leq n,  \\
       	& 0 \leq y_i \leq 1,  \; \mbox{ for all } 1\leq i \leq n,  \\
\end{array}
 \end{equation}
where $e\in\mathbb{R}^{n}$ denotes the vector of ones. We note that the last $n$ constraints denote
 a simple box in the auxiliary variable vector $y\in\mathbb{R}^{n}$.  A more difficult  reformulation
  substitutes the simple box by a set of binary constraints given by:  either $y_i = 0$ or $y_i =1$ for all $i$.
In that case, the problem is an integer programming problem (much harder to solve) for which there are several
 algorithmic ideas already developed; see, e.g., \cite{Bertsimas09, Bienstock, Cesarone2013, DiLorenzo, ZhengSunLi}. In here, we will focus on
 the continuous formulation (\ref{relaxed}), that will play a key role in our algorithmic proposal. For additional theoretical properties that
 include the equivalence between the original version (\ref{origprb}) and the continuous
   relaxed version (\ref{relaxed}), see \cite{Burdakov,Kanzow,Krul1,Krul2}. \\

As a consequence of the so-called Hadamard constraint ($x\circ y =0$, i.e., $x_i y_i = 0$ for all $i$),
the formulation (\ref{relaxed})  is a nonconvex
 problem, even when the original cardinality constrained problem (except for the cardinality constraint of course)
  was convex. Thus, one can in general not expect to obtain global minima. But if one is for example
interested in obtaining local solutions or good starting points for a global method, this continuous formulation (\ref{relaxed}) can be useful. \\

In this work, we will pay special attention to those problems for which the set
$\Omega$ is the intersection of a finite collection of convex sets, in such a way that it is very easy to project
 onto each one of them. In that case, the main idea is to take advantage of the fact that two of the
  constraints in (\ref{relaxed}), namely $e^{\top} y \geq n-\alpha$ and $0 \leq y_i \leq 1$ for all $i$,
   are also ``easy-to-project" convex sets, and so an alternating projection scheme
   can be conveniently applied to project onto the intersection of all the involved constraints in (\ref{relaxed}), except
    for the Hadamard constraint. For solving the continuous formulation (\ref{relaxed}) we can then use a suitable  low-cost
    convex constrained scheme, such as gradient-type methods in which  the objective function includes $f(x)$ plus a suitable
    penalization term that guarantees that the Hadamard constraint is also satisfied at the solution. In Section \ref{penaltysec},
    we will describe and analyze a  general penalty method to satisfy the Hadamard constraint that appears in the relaxed formulation (\ref{relaxed}).
    In Section \ref{spgdyks}, we will describe a suitable alternating projection scheme as well as a suitable low-cost gradient-type projection
     method  that can be combined with the penalty method of Section \ref{penaltysec}.  Concerning some specific applications, in Section \ref{portfsec},
   we will consider in detail the standard mean-variance limited diversified portfolio selection problem
    (see e.g., \cite{Cesarone2009, Cesarone2013, ChangEtAl, DiLorenzo, Fastr, Gao}). In Section \ref{compexp},  we will present a numerical study to
    illustrate the computational performance of the proposed scheme on a variety of data sets involving real-world capital
 market indices from major stock markets. For each considered data set, we will focus our attention on the efficient frontier produced for different
  values of the limited number of allowed assets.  In Section \ref{conclu}, we will present  some final comments and perspectives. \\

   \section{A penalization strategy for the Hadamard constraint} \label{penaltysec}

   Let us consider again the continuous  formulation (\ref{relaxed}), and let us focus our attention on the Hadamard
    constraint $x\circ y  = 0$ (i.e., $x_i y_i = 0$ for all $i$). This particular constraint,
    is the only one that does not define a convex set. The others define convex sets in which it is easy
     to project, as discussed in the previous section. To see that the set of vectors $(x,y) \in \R^{2n}$
     such that $x\circ y  = 0$ do not form a convex set, it is enough to consider the two 2-dimensional
    pairs: $(x_1,y_1) = (1,0,0,1)$ and $(x_2,y_2) = (0,1,1,0)$. Both pairs are clearly in
      that set, but the convex combination: $\frac{1}{2}(x_1,y_1) + \frac{1}{2}(x_2,y_2) = \frac{1}{2}e$, which
   is not in that set. \\

  A classical and straightforward approach to force the Hadamard condition at the solution, while keeping the
   feasible set of our problem as  the intersection of a finite collection  of easy convex sets,  is to add a penalization term $\tau h(x,y)$ to the objective function and consider instead the following formulation:
    \begin{equation} \label{relaxed2}
  	 \begin{array}{ll}
\displaystyle	\min_{x,y}   & f(x) + \tau h(x,y) \\
	\mbox{subject to:} &  x\in\Omega,  \\
       	& e^{\top} y \geq n-\alpha,  \\
       	& 0 \leq y_i \leq 1,  \; \mbox{ for all } 1\leq i \leq n,  \\
\end{array}
 \end{equation}
  where $\tau >0$  is a penalization parameter that needs to be properly chosen, and the function $h:\mathbb{R}^{2n}\to \mathbb{R}$ is  continuously differentiable and chosen to satisfy the following two properties: $h(x,y)\geq 0$ for all feasible vectors $x$ and $y$, and
  $h(x,y) = 0$ if and only if $x\circ y  = 0$. Clearly, the function $h(x,y)$ is crucial and should be
   conveniently chosen depending on the considered application. \\

   Applying now a penalty scheme, problem  (\ref{relaxed2}) can be reduced to a sequence of convex constrained problems of the following form:
  \begin{equation}\label{pnlty}
  \displaystyle	\min_{x,y} \;\;\;\; f(x) + \tau_k h(x,y), \;\;\;\; \mbox{ subject to }\;\;\;\;
   (x,y) \in\widehat{\Omega},
  \end{equation}
  where $\tau_k >0$ is the penalty parameter that  increases at every $k$ to penalize the
  Hadamard-constraint violation, and the closed convex set $\widehat{\Omega}$ is given by
 \[ \widehat{\Omega} = \{(x,y)\in \R^{2n} : x\in\Omega, \;\; e^{\top} y \geq n-\alpha, \;\;  0 \leq y_i\leq 1, \; i=1,\ldots,n\}. \]

   Under some mild assumptions and some specific choice
  of  the sequence $\{\tau_k\}$, it can be established that the sequence of solutions of  problem (\ref{pnlty}) converges to a solution of
   (\ref{relaxed}); see, e.g.,  \cite{Fiacco68} and   \cite[Secc. 12.1]{Luenb84}.
 Let us assume that problem (\ref{relaxed}) attains global minimizers.
 Since $f$ is a continuous function, it is enough to assume that one of the closed and convex sets involved in
  the definition of $\Omega$  in (\ref{relaxed}) is bounded.  In here, for the sake of completeness,
  we summarize the convergence properties of the proposed  penalty scheme (\ref{pnlty}).
   \begin{theorem} \label{penaltyt}
   If for all $k$, $\tau_{k+1} > \tau_k>0$  and $(x_k,y_k)$ is a global solution of (\ref{pnlty}), then
  \begin{eqnarray*}
  f(x_k) + \tau_k h(x_k,y_k) & \leq & f(x_{k+1}) + \tau_{k+1} h(x_{k+1},y_{k+1}) \\
  h(x_{k+1},y_{k+1}) & \leq &  h(x_k,y_k) \\
  	f(x_k) & \leq & f(x_{k+1}) \;.
  \end{eqnarray*} 	
 Moreover, if $\bar{x}$ is a global solution of  problem (\ref{relaxed}), then for all $k$
 $$  f(x_k) \;\leq\;  f(x_k) + \tau_k h(x_k,y_k) \;\leq\; f(\bar{x}) \;.  $$  	
 Finally,  if $\tau_k \rightarrow \infty$ and  $\{(x_k,y_k)\}$ is the sequence of global minimizers
 obtained by solving (\ref{pnlty}) then any limit point of  $\{(x_k,y_k)\}$ is a global minimizer of
   (\ref{relaxed}).
    \end{theorem}
 \begin{remark} \label{taubdd}
In the proof of the last statement of Theorem \ref{penaltyt} (see, e.g., \cite[Secc. 12.1]{Luenb84}), the requirement of $\tau_k \rightarrow \infty$
is used only to guarantee that the term $h(x_k,y_k)\rightarrow 0$ when $k\rightarrow \infty$, i.e., to
 guarantee that $x_k\circ y_k \rightarrow 0$. In order to guarantee the convergence result, what is important  is that the
  Hadamard product itself goes to zero even if $0< \tau_k < \infty$ for all $k$. This fact will play a key role in our numerical study (Section 5).  	
 \end{remark}   	

We would like to close this section with a pertinent result (\cite[Theorem 4]{Krul1}) that establishes a
one-to-one correspondence between minimizers of problems (\ref{origprb}) and  (\ref{relaxed}), whenever
the obtained solution $\bar{x}$ satisfies the cardinality constraint with equality,  i.e.,  $\|\bar{x}\|_0=\alpha$.
\begin{theorem}\label{theo: min_bin}
	Let $(\bar{x},\bar{y})$ be a local minimizer of the relaxed problem (\ref{relaxed}). Then
 $\|\bar{x}\|_0=\alpha$ if and only if $\bar{y}$ is unique, that is, if there exist exactly one $\bar{y}$
 such that $(\bar{x},\bar{y})$ is a local minimizer of (\ref{relaxed}). In this case, the components of
$\bar{y}$ are binary (i.e., $\bar{y}_i =0$ or $\bar{y}_i=1$ for all $1\leq i\leq n$) and $\bar{x}$ is a local minimizer of (\ref{origprb}).
\end{theorem}

  \section{Dykstra's method and the SPG method} \label{spgdyks}

 For every $k$, a low-cost projected gradient method can be used to solve the optimization problem (\ref{pnlty}).
 Notice that $\widehat{\Omega}$ is the intersection of finitely many ``easy" convex sets. A convenient tool for finding the required
   projections onto $\widehat{\Omega}$ is Dykstra's alternating projection algorithm  \cite{BoyleDykstra1986},
 that will be now described in a general setting in $\R^n$.
 Roughly speaking, Dykstra's algorithm projects in a clever way onto the easy
 convex sets individually to  complete a cycle which is repeated iteratively, and as any other iterative
  method it can be stopped prematurely. \\

  For a given $\tilde{x}\in \R^n$ and finitely many closed convex sets, say  $\Omega_1,\dots,\Omega_p$ in $\R^n$, we consider the best approximation 
  problem: find the closest point to $\tilde{x}$ in
  $\widehat{\Omega}=\cap_{i=1}^{p}\Omega_i \ne \emptyset$,  which  can be stated as an optimization problem as
   follows:
  \begin{equation} \label{Project}
  	\mbox{minimize}\,\|\tilde{x}-x\| \quad\mbox{ subject to}\quad x\in \widehat{\Omega},
  \end{equation}
  where, for any $z\in \R^n$, $\|z\|^2 = \left\langle z,z\right\rangle$. The unique solution
  $x^{*}$ of problem (\ref{Project}) is called the projection of $\tilde{x}$ onto $\widehat{\Omega}$ and is denoted as $P_{\widehat{\Omega}}(\tilde{x})$.\\

  In Dykstra's method it is assumed that the
  projections onto each of the individual sets  $\Omega_i$ are relatively simple to compute, e.g., boxes, spheres, subspaces, half-spaces, hyperplanes, among others.
  The algorithm has been adapted and used for solving a huge amount of different applications.
  For a review on Dykstra's  method, its properties and applications, as well as many other alternating projection schemes; see, e.g., \cite{Deutsch, EyR}.\\

  Dykstra's algorithm solves (\ref{Project}) by generating two sequences: the iterates $\{ x_{\ell}^{i} \}$ and the increments $\{I_{\ell}^{i}\}$. These sequences are defined by the following recursive formulae:
  \begin{equation} \label{alg-Dyks}
  	\begin{array}{lclll}
  		x^{0}_\ell & = & x_{\ell-1}^p & \\ [2mm]
  		x_{\ell}^i & = & P_{\Omega_i}(x_{\ell}^{i-1}-I_{\ell-1}^i) & \; &
  		i=1,2, \ldots, p, \\ [2mm]
  		I_{\ell}^i & = & x_{\ell}^i-(x_{\ell}^{i-1}-I_{\ell-1}^i) & \; &
  		i=1,2, \ldots, p,
  	\end{array}
  \end{equation}
  for $\ell \in\Z^+$ with initial values $x_0^p=\tilde{x}$ and $I_0^i=0$  for $i=1,2, \dots, p$. \\

  Notice that the increment $I_{\ell-1}^{i}$ associated with $\Omega_i$ in the
  previous cycle is always subtracted before projecting onto
  $\Omega_i$. The  sequence of increments  play a fundamental role in the convergence
  of the sequence $\{ x^{i}_{\ell} \}$ to the  unique optimal solution  $x^{*}=P_{\widehat{\Omega}}(\tilde{x})$
  of problem (\ref{Project}).
    Notice also that, for the sake of simplicity in our presentation, the projecting
  cyclic control index $i(\ell)$ used in (\ref{alg-Dyks}) is the most common one: $i(\ell) = \ell \mbox{
  	mod }p + 1$, for all $\ell\geq 0$. However, more advanced control
  indices can also be used, as long as they satisfy some minimal theoretical requirements; see, e.g., \cite{EyR}). \\

   Boyle and Dykstra \cite{BoyleDykstra1986} established the key convergence theorem associated
  with  algorithm (\ref{alg-Dyks}).
  \begin{theorem} \label{BoDyk}
  	Let $\Omega_1,\dots ,\Omega_p$ be closed and convex sets of $\R^n$
  	such that $\widehat{\Omega}=\cap_{i=1}^p \Omega_i \neq \emptyset$. For any
  	$i=1, 2, \ldots,p$ and any $\tilde{x} \in \R^n$, the sequence $\{
  	x_{\ell}^i \}$ generated by (\ref{alg-Dyks}) converges to
  	$x^{*}=P_{\widehat{\Omega}}(\tilde{x})$ (i.e., $\|x_{\ell}^{i}-x^{*}\| \rightarrow 0$
  	as $\ell \rightarrow \infty $).
  \end{theorem}
  Concerning the rate of convergence, it is well-known that  Dykstra's algorithm exhibits a  linear rate of
   convergence in the polyhedral case (\cite{Deutsch, EyR}), which is the case in all problems considered here, see Section 5.
   Finally, the stopping criterion associated with Dykstra's algorithm is a delicate issue. A discussion about
    this topic and the development of some robust stopping criteria are fully described in \cite{BR2005}.
   Based on that, in here we will stop the iterations when
 \begin{equation} \label{stopdyk}
 \sum^p_{i=1} \| I_{\ell-1}^i-I_\ell^i \|^{2}\leq \varepsilon,
 \end{equation}
  where $\varepsilon > 0$ is a small given tolerance. \\

 Since the gradient $\nabla f(x,y)$ of $f(x,y) = f(x) + \tau h(x,y)$ is available for each fixed $\tau >0$,
  then Projected Gradient (PG) methods provide an interesting low-cost option for solving (\ref{pnlty}). They are
  simple and  easy to code, and avoid the need for matrix factorizations (no Hessian matrix is used).
There have been many different variations of the early PG methods. They all have the common property of
 maintaining feasibility of the iterates by frequently projecting trial steps on the feasible convex set.
  In particular, a well-established and effective scheme is the so-called Spectral Projected Gradient (SPG) method; see Birgin et al.
  \cite{bmr1,bmr2,bmr3, bmr4}).  \\

 The SPG algorithm starts with $(x_0,y_0) \in \R^{2n}$, and  moves at every iteration $j$ along the internal
 projected gradient  direction $d_j=P_{\widehat{\Omega}}((x_j,y_j) -\alpha_j \nabla f(x_j,y_j)) - (x_j,y_j)$,
 where $d_j \in \R^{2n}$ and $\alpha_j$ is the well-known spectral choice of step length (see \cite{bmr4}):
 \[ \alpha_j = \frac{\langle s_{j-1}, s_{j-1}\rangle}{\langle s_{j-1},
 	(\nabla f(x_j,y_j) - \nabla f(x_{j-1},y_{j-1})) \rangle}, \]
 and $s_{j-1}= (x_j,y_j) -  (x_{j-1},y_{j-1})$.
 In the case of rejection of the first trial point, $(x_j,y_j)+d_j$, the
 next ones are computed along the same direction, i.e., $(x_{+},y_{+})=(x_j,y_j) + \lambda d_j$, using a
  nonmonotone  line search to choose  $0<\lambda\leq 1$ such that the following condition holds
 \[
 f(x_{+},y_{+}) \leq \max_{0\leq l \leq \mbox{ min } \{j,M-1\}}
 f(x_{k-l},y_{k-l}) + \gamma \lambda \langle d_j, \nabla f(x_j,y_j) \rangle, \]
 where $M \geq 1$ is a given integer and $\gamma$ is a small
 positive number. Therefore,  the projection onto $\widehat{\Omega}$  must
 be performed only once per iteration. More details can be found in
 \cite{bmr1} and \cite{bmr2}. In practice $\gamma = 10^{-4}$ and  a typical value for the
 nonmonotone parameter is  $M=10$, but the performance of the method may vary for variations of
 this parameter and a fine tuning may be adequate for specific  applications.\\

  Another key feature of the SPG method is to accept the initial spectral step-length as often as possible
  while ensuring global convergence. For this reason, the SPG method employs a non-monotone
   line search that does not impose functional decrease at every
   iteration.  The global convergence of  the SPG method combined with Dykstra's algorithm to
   obtain the required projection per iteration can be found in   \cite[Section 3]{bmr3}. \\

\section{Cardinality constrained optimal portfolio problem} \label{portfsec}

Let the vector $v\in \R^n$ and the symmetric and positive semi-definite matrix $Q\equiv[\sigma_{ij}]_{i,j=1,\ldots,n}\in \R^{n\times n}$ be
the given mean return vector and  variance-covariance matrix of the $n$ risky available assets, respectively. The entry
  $\sigma_{ij}$ in $Q$ is the covariance between assets $i$ and $j$ for $i,j=1,\ldots,n$ , $\sigma_{ii}=\sigma_i^2$ and $\sigma_{ij}=\sigma_{ji}$.
  As a consequence of the pioneering work of Markowitz \cite{Markowitz}, the   mean-variance
 portfolio selection problem can be formulated as (\ref{origprb}), where the  objective function is given by
 \begin{equation} \label{risk}
 f(x) = \frac{1}{2} x^{\top}Qx,
  \end{equation}
  and the convex set
 \[ \Omega = \{x\in \R^n : v^{\top} x \geq \rho, \;\; e^{\top} x=1, \;\;  0 \leq x_i\leq up_i, \; i=1,\ldots,n\}, \]
 representing the constraints of minimum expected return level $\rho$, budget constraint ($\sum_{i=1}^n x_i = 1$
 means that all available wealth will be invested), and lower ($x\geq0$ excludes short sale) and upper bounds for each $x_i$, respectively. Notice
  that the minimization of $f(x)$, involving the given covariance matrix $Q$, accounts for the minimization of the
  variance,
  while the return is expected to be at least $\rho$. Notice also that, as
   previously discussed,  in this case the set $\Omega$ is the intersection of three easy convex sets: a half-space, a hyperplane, and a box.
 The additional constraint in (\ref{origprb}), $\|x\|_0 \leq \alpha$ for $0<\alpha<n$, plays a key role here, and
  indicates that among the $n$ risky available options, we can only invest in at most $\alpha$ assets
  (cardinality constraint). The solution vector $x$ denotes an investment portfolio and each $x_i$ represents the fraction held of each
  asset $i$.  It should be mentioned that other inequality and/or
  equality constraints can be added to the problem, as they represent additional real-life constraints; e.g., transaction costs \cite{Bertsimas, KKR}. \\

  Now, as discussed above, our main idea is to consider the continuous formulation (\ref{relaxed}) instead of
  the optimization problem (\ref{origprb}). For the portfolio selection problem  we would end up
   with the following problem that involves the auxiliary vector $y$:
   \begin{equation} \label{portfol}
  	\begin{array}{ll}
  		\displaystyle\min_{x,y}   & \frac{1}{2} x^{\top}Qx \\ [1mm]
  		\mbox{subject to:} & v^{\top} x \geq \rho,  \\ [1mm]
  		& e^{\top} x = 1,  \\ [1mm]
  		& 0 \leq x_i\leq up_i,  \; \mbox{ for all } 1\leq i \leq n,  \\ [1mm]
  		& e^{\top} y \geq n-\alpha,  \\ [1mm]
  		& x\circ y = 0,   \\ [1mm]
  		& 0 \leq y_i \leq 1,  \; \mbox{ for all } 1\leq i \leq n,
  	\end{array}
  \end{equation}
 where the upper bound vector $up\in \R^n$  and $\rho>0$ are given. Note that the vector
  $y$ appears only in the last 3 constraints, and the vector $x$ appears in the first three constraints
   but also in the (non-convex) Hadamard constraint: $x\circ y  = 0$.  \\

    As discussed in Section 2, the best option to force the Hadamard  condition at the
        solution while keeping the feasible set of our problem as  the intersection of a finite collection
   of easy convex sets,  is to add the term $\tau h(x,y)$ to the objective function, where our convenient choice is $h(x,y)= x^{\top}y$:
   \begin{equation} \label{minmaxhad}
  	f(x,y) = \frac{1}{2} x^{\top}Qx + \tau x^{\top}y,
  \end{equation}
  where $\tau >0$
   is a penalization parameter that needs to be properly chosen as described in Section 2.  Since the vectors $x$ and $y$ will be forced
    by the alternating projection scheme to have all their entries greater than or equal to zero, then $h(x,y)= x^{\top}y\geq 0$ for any
   feasible pair $(x,y)$,  and forcing $\tau x^{\top}y=0$ is equivalent to forcing the Hadamard condition: $x_i y_i = 0$ for all $i$.
   Notice that, setting $\tau=0$ for solving  (\ref{portfol}) with $f(x,y)$ given by (\ref{minmaxhad})  minimizes the risk,
   independently of the Hadamard condition. On the other hand, if $\tau>0$ is sufficiently large as compared to the size of $Q$ then the
    term $x^{\top}y$ must be zero at the solution. Hence,  choosing $\tau>0$ represents an explicit trade-off between the risk and
    the Hadamard condition.  \\

 Our algorithmic proposal consists in solving a sequence of penalized problems, as described in Section 2, using the SPG scheme and Dykstra's
  alternating projection method (that from now on will be denoted as the SPG-Dykstra method) to solve problem (\ref{portfol}), without the complementarity constraint $x\circ y=0$, and using the objective function given by (\ref{minmaxhad}).  That is, for a sequence of increasing penalty terms
  $\tau_k>0$, we will solve the following problems
    \begin{equation} \label{portfol2}
	\begin{array}{ll}
	\displaystyle	\min_{x,y}   & \frac{1}{2} x^{\top}Qx + \tau_k x^{\top}y\\ [1mm]
		\mbox{subject to:} & v^{\top} x \geq \rho,  \\ [1mm]
		& e^{\top} x = 1,  \\ [1mm]
		& 0 \leq x_i\leq up_i,  \; \mbox{ for all } 1\leq i \leq n,  \\ [1mm]
		& e^{\top} y \geq n-\alpha,  \\ [1mm]
		& 0 \leq y_i \leq 1,  \; \mbox{ for all } 1\leq i \leq n.
	\end{array}
\end{equation}
Since the function $h(x,y)= x^{\top}y$ satisfies the properties mentioned in section 2, if we choose the sequence of parameters $\{\tau_k\}$ such
that $h(x_k,y_k)$ goes to zero when $k$ goes to infinity, then Theorem \ref{penaltyt} guarantees the convergence of the proposed scheme.

   Before showing some computational results in our next section, let us recall that the gradient and the Hessian of the objective function $f$
   at every pair $(x,y)$ are given by
\begin{equation*}
	\nabla f(x,y) = \left(\begin{array}{c}
		Qx + \tau_k y	\\
		\tau x\end{array}\right)\;\; \mbox{ and }\;\;
	\nabla^2 f(x,y) = \left(\begin{array}{cc}
		Q & \tau_k I \\
		\tau_k I & 0 	\end{array}\right).
\end{equation*}
Notice that, for any $\tau_k >0$,  $\nabla^2 f(x,y)$ is symmetric and indefinite.

\section{Computational results} \label{compexp}

To add understanding and illustrate the advantages of our proposed combined scheme, we present the results of some numerical experiments on an
 academic simple problem ($n=6$), and also on some data sets involving real-world capital market indices from major stock markets.
 All the experiments were performed using Matlab R2022 with double precision on an Intel$^{\circledR}$ Quad-Core i7-1165G7
 at 4.70 GHz with 16GB of RAM memory,  using Windows 10 Pro with 64 Bits.
 \\

The algorithm we use in this section was indicated in the previous sections and now, for completeness, we describe it in detail.

 \noindent{\bf { {Algorithm Penalty-SPG-Dykstra (PSPGD)}}. }
 \begin{itemize}
 \item[S0]: { {Given $\tau_{-1} >0$}}, set $x_{-1}=(1/n)e$,  $y_{-1} = 0$, and $k=0$.

 \item[S1]: Compute { {$ \tau_k > \tau_{k-1}$}}
 \item [S2]: Set $ x_{k,0} = x_{k-1} \mbox{ and } y_{k,0} = y_{k-1}$,  and { {from $(x_{k,0},y_{k,0})$ apply the SPG-Dykstra method to
 (\ref{portfol2}),  until}}
 $$\|P_{\widehat{\Omega}}((x_{k,m_{k}},y_{k,m_{k}}) - \nabla f(x_{k,m_{k}},y_{k,m_{k}})) - (x_{k,m_{k}},y_{k,m_{k}})\|_2\leq { tol_1} $$
 is satisfied at some  iteration $m_{k}\geq 1$. Set $ x_{k} = x_{k,m_{k}} $ and $ y_{k} = y_{k,m_{k}}$.
 \item[S3]: If
 $$ x_{k}^{\top}y_{k} \leq { tol_2} \;\mbox{ and }\;  |f(x_{k})-f(x_{k-1})| \leq  tol_2 $$ then stop. Otherwise,
 set $ k= k+1 $ and return to S1.
 \end{itemize}

 For our experiments, we set $tol_1=10^{-6}$ and  $tol_2=10^{-8}$.
  We note that at any iteration $k\geq 1$,  Step S2 of Algorithm PSPGD starts from $(x_{k-1},y_{k-1})$, which is the previous
solution of  (\ref{portfol2}),  obtained  using $\tau_{k-1}$.
 We also note that to stop  the SPG-Dykstra iterations we monitor the value of $\|P_{\widehat{\Omega}}((x_k,y_k) - \nabla f(x_k,y_k)) - (x_k,y_k)\|_2$,
  which is denoted as the pgnorm at iteration $k$ in the tables below.
  It is worth recalling that  if $\|P_{\widehat{\Omega}}((x,y) - \nabla f(x,y)) - (x,y)\|_2=0$, then
   $(x,y) \in \widehat{\Omega}$ is stationary for  problem (\ref{portfol2}); see, e.g.,  \cite{bmr1,bmr3}.
   Concerning the nonmonotone line search strategy used by the SPG method, we set $\gamma = 10^{-4}$ and  $M=10$.
    Each SPG iteration uses Dykstra's altrnating projection scheme to obtain the required projection onto $\widehat{\Omega}$, and this internal
    iterative process is stopped when (\ref{stopdyk}) is satisfied with $\varepsilon = 10^{-8}$. \\

 To explore the behavior of Algorithm PSPGD, we will vary  the minimum expected return parameter $\rho>0$  and the cardinality constraint
  positive integer $1\leq \alpha <n$.  In all cases, we set the upper bound vector $up=e$, where $e$ is the vector of ones.
  Of course, for certain combinations of all those parameters the problem might be
   infeasible.  We will discuss possible choices of these parameters to guarantee that the feasible region
    of problem (\ref{portfol2}) is not empty. \\

 To keep a balanced trade-off between the risk and
    the Hadamard condition, it is convenient to choose the initial parameter { {$\tau_{-1}>0$}} of the same order of
    magnitude  of the largest eigenvalue of $Q$. For that, we proceed as follows: set $z= Qe$ and
    ${ {\tau_{-1}}} = z^{\top}Qz/(z^{\top}z)$, i.e.,
   a Rayleigh-quotient of $Q$ with a suitable vector $z$, which produces a good estimate of $\lambda_{\max}(Q)$.  This choice worked well for the 
   vast majority of the test examples.  According to  Remark \ref{taubdd}, to observe convergence, we need to drive the inner product
    $x_k^{\top}y_k$ down to zero.  For that   we increase the penalization parameter as follows:
\begin{equation} \label{tauk}
	\tau_{k+1}=\delta_{k+1}\tau_k\;\; \mbox{ where }\;\; \delta_{k+1}=\delta_k+\frac{(n-\alpha)\rho}{n}\dfrac{|v^{\top}x_{k+1}|}{\sqrt{x_{k+1}^{\top}Qx_{k+1}}}
  \;\; \mbox{ and } \;\; \delta_{-1}=1.
\end{equation}
We note that in practice this formula increases the penalty parameter  in a controlled way taking into account the ratio between the absolute value of
the current return $|v^{\top}x_{k+1}|$ and the current risk $\sqrt{x_{k+1}^{\top}Qx_{k+1}}$.  In all the reported experiments, the controlled
sequence $\{\tau_k\}$  given by (\ref{tauk}) was enough to guarantee that the Hadamard product goes down to zero. \\

 Concerning the choice of the expected return, based on \cite{Cesarone2009,ZhengSunLi}, in order to consider feasible
  problems we study the behavior of our combined scheme in an interval $[\rho_{\min}, \rho_{\max}]$ of possible values of the parameter $\rho$, which
   is obtained as follows.  Let $\rho_{\min}=v^{\top}x_{\min}$ and $\rho_{\max}=v^{\top}x_{\max}$ where
\begin{align*}
	x_{\min}	=\arg \min_{x}&\; \dfrac{1}{2}x^{\top}Qx+\tau x^{\top}y\\
	\mbox{ subject to:}&\;\; e^{\top}x=1,\\
	&\;\; 0 \leq x_i\leq up_i,  \; \mbox{ for all } 1\leq i \leq n,  \\
	&\;\; e^{\top} y \geq n-\alpha,  \\
	&\;\; 0 \leq y_i \leq 1,  \; \mbox{ for all } 1\leq i \leq n,
\end{align*}
and
\begin{align*}
	x_{\max}	=\arg \max_{x} &\; v^{\top}x-\tau x^{\top}y \\
	\mbox{ subject to:}&\;\; e^{\top}x=1, \\
	&\;\; 0 \leq x_i\leq up_i,  \; \mbox{ for all } 1\leq i \leq n,  \\
	&\;\; e^{\top} y \geq n-\alpha,  \\
	&\;\; 0 \leq y_i \leq 1,  \; \mbox{ for all } 1\leq i \leq n.
\end{align*}
 These two auxiliary optimization problems are solved in advance, only once for each considered problem, using in turn the proposed
 Algorithm PSPGD. For that, we fix the same parameters and we start from the same initial values indicated above.
 Once the interval $[\rho_{\min}, \rho_{\max}]$ has been obtained, to choose a suitable return $\rho$ we can proceed as follows.
 For a fixed $0<\tilde{\epsilon}<1$, if $\rho_{\min}+\tilde{\epsilon}(\rho_{\max}-\rho_{\min})\geq 0$ we set $\rho=\rho_{\min}+\tilde{\epsilon}(\rho_{\max}-\rho_{\min})$,
 else if  $|\rho|\leq v_{\max}$ we set $\rho=\tilde{\epsilon}|\rho|$, otherwise we set $\rho=\tilde{\epsilon} v_{\max}$. In here,
 $v_{\min}=\min\{v_1,\ldots,v_n\}$ and $v_{\max}=\max\{v_1,\ldots,v_n\}$. 	\\

For our first data set we consider a simple portfolio problem with $n=6$ available assets, denoted as {\it Simple-case} for which the
 mean return vector $v$ and the covariance matrix $Q$ are given by:
$$ 	v=(0.021\;\; 0.04\;\; -0.034\;\; -0.028\;\; -0.005\;\; 0.006)^{\top}, $$
 \[ Q =  \left[ \begin{array}{rrrrrr}
 	0.038 & 0.020 & 0.017 &  0.014 &  0.019 &  0.017 \\
 	0.020 & 0.043 & 0.015 &  0.013 & 0.021 & 0.014 \\
   	0.017 & 0.015 & 0.034 & 0.011 & 0.014 & 0.014 \\
   	0.014 & 0.013 & 0.011 & 0.044 & 0.014 & 0.011 \\
   0.019 & 0.021 &  0.014 & 0.014 &  0.040 &  0.014 \\
   0.017 & 0.014 & 0.014 & 0.011 & 0.014 & 0.046 \end{array} \right]. \]
We note that $Q$ is symmetric and positive definite ($\lambda_{\min}(Q) =1.79 \times 10^{-2}$ and
	 $\lambda_{\max}(Q)=1.17\times 10^{-1}$). Notice that the assets three, four, and five have negative average returns.
 The purpose of this simple example is to demonstrate properties of the problem and the proposed algorithm in an easy-to follow fashion.
 For the other data sets, involving real-world capital market indices, we consider  some larger problems obtained from Beasley's OR Library
(http://people.brunel.ac.uk/$\!_{^{\sim}}\!$mastjjb/jeb/info.html), built from weakly price data from March 1992
 to September 1997, and that we will denote as Port1 (Hang Seng index with $n=31$), Port2 (DAX index with $n=85$),
 Port3 (FTSE 100 index with $n=89$),  Port4 (S\&P 100 index with $n=98$), Port5 (Nikkei index with $n=225$), and
  Port 6 ($n=600$, former by assets from NY Stock Exchange, weekly prices from July 2001 to July 2018, \cite{JKMP}); see also \cite{Beasley,ChangEtAl}. \\

  The key properties, to be discussed and illustrated in the rest of this section, are the influence of the cardinality constraint to the feasible
   set in the risk-return plane, the efficient frontier, and the quality of the solution obtained by Algorithm PSPGD. The feasible
   set is usually represented in the risk-return plane, presenting all possible combinations of assets that satisfy the constraints. In general the
    feasible set for the classical problem without cardinality constraint has the so-called bullet shape.  The efficient frontier is the set of
     optimal portfolios that offer the highest expected return for a defined level of risk or the lowest risk for a given level of expected return.
      Clearly, in the risk-return plane, the efficient frontier is the upper limit of the feasible set. \\

 Introducing the cardinality constraints might complicate the feasible set in the sense that the set is shrinking as we will now show.
 Starting with the feasible interval for the expected return we report in Table~\ref{Table:rho}, $\rho_{\max}\leq v_{\max}$ and
$\rho_{\min} \geq  v_{\min}$,  for $\alpha=5$ and for all the considered data sets.

\begin{table}[H]
	\begin{center}
		\begin{tabular}{cccccc}
			Problem &   $n$ &   $v_{\min}$ &$v_{\max}$ &$\rho_{\min}$ & $\rho_{\max}$ \\ \hline \hline
			Simple case & 6  & -0.0340& 0.0400 & -0.0238&0.0373 \\ \hline
			Port1 & 31 &   5.64e-4    & 0.0435& 0.0130&0.0435 \\ \hline
			Port2 & 85 &    -0.0160   & 0.0392&0.0099&0.0342\\\hline
			Port3 & 89 &  -0.0045   & 0.0328&0.0102&0.0268\\\hline
			Port4 & 98 &  -0.0079   & 0.0368&0.0077&0.0271\\\hline
			Port5 & 225 &  -0.0340   & 0.0159&-0.0060&-0.0060\\\hline
			Port6 & 600 &  -0.0593   &0.0364&0.0013& 0.0013\\\hline
		\end{tabular}
	\end{center}
	\caption{Return value with $ \alpha = 5 $ for all data sets.}\label{Table:rho}
\end{table}

\begin{figure}[H]
	\begin{center}
		\includegraphics[width=18cm, height=9.5cm]{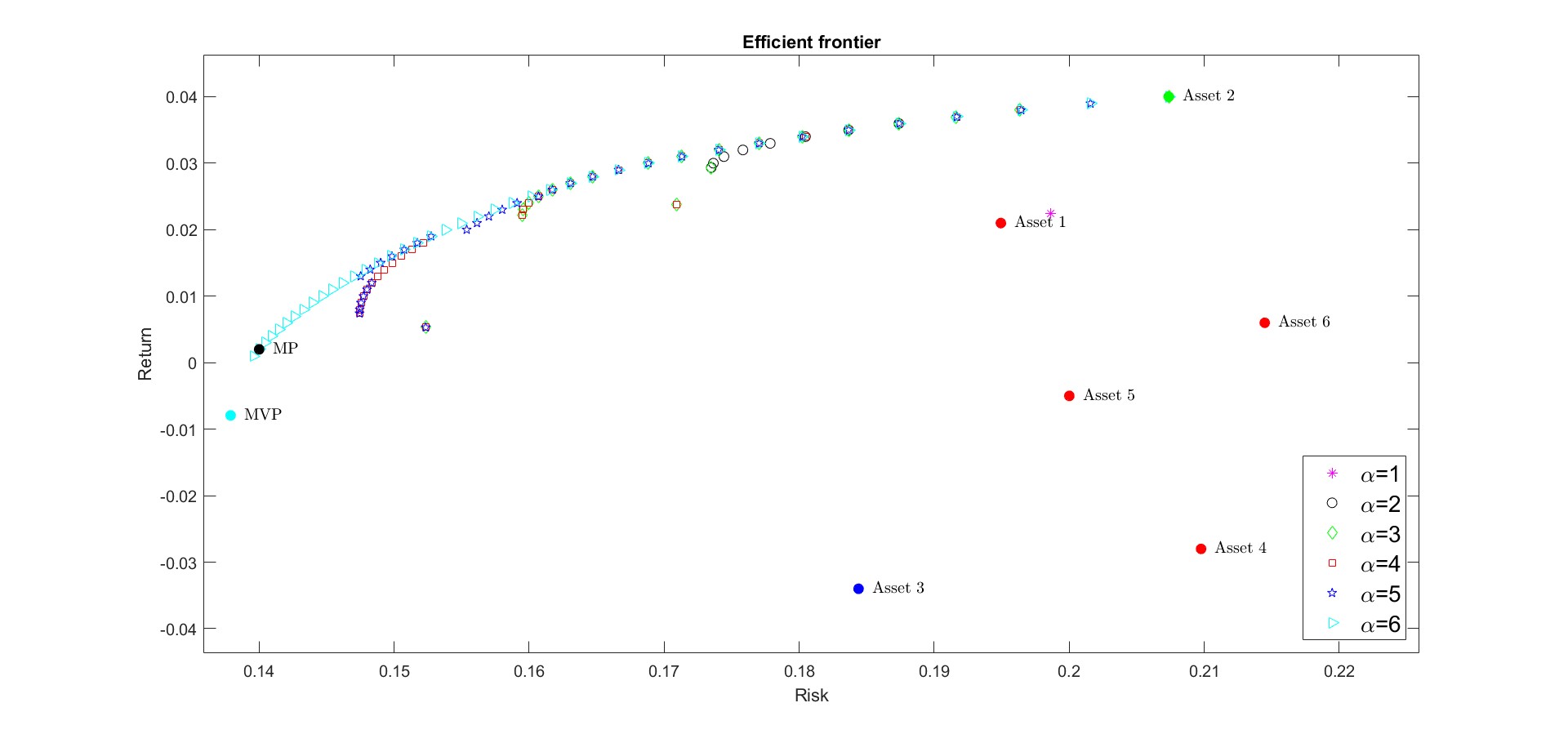}
	\end{center}
	\caption{Risk versus return, using Algorithm PSPGD for the Simple-case problem.}
\label{fig14}
\end{figure}

Let us now take a closer look at the Simple-case.
If we solve the original Markowitz problem \cite{Markowitz} - the minimal variance portfolio,  (i.e.,  $\displaystyle\min_{x}  \frac{1}{2} x^{\top}Qx$
subject to $e^{\top} x = 1$) for the Simple-case problem  we obtain
$$\bar{x}=(0.0961,0.1168,0.2625,0.2140,0.1429,0.1677)^{\top}, $$
 risk $\sqrt{\bar{x}^{\top}Q\bar{x}}=0.1379$, and expected return $v^{\top}\bar{x}=-0.0079$. Solving the same problem with the additional constraint
  $x\geq0$ we get the same solution. Thus, the minimal variance portfolio is the same as the minimal variance portfolio without short sale.
In Figure~\ref{fig14}, we present for the Simple-case problem, the return and risk for  all 6 assets, the minimal variance portfolio, denoted by MVP, the classical
 Markowitz portfolio without short sale and the expected return  constraint  $v^{\top}x\geq \rho=0.002$, denoted by MP, as well as the efficient frontier for
 different values of the cardinality constraint $ \alpha$.   Clearly for $ \alpha = 6$, i.e., without cardinality constraint, we get a classical
 convex efficient frontier while for smaller  $\alpha$ values the curves are deformed. \\

\begin{figure}[H]
	\begin{center}\hspace{-2cm}
		\includegraphics[width=19cm]{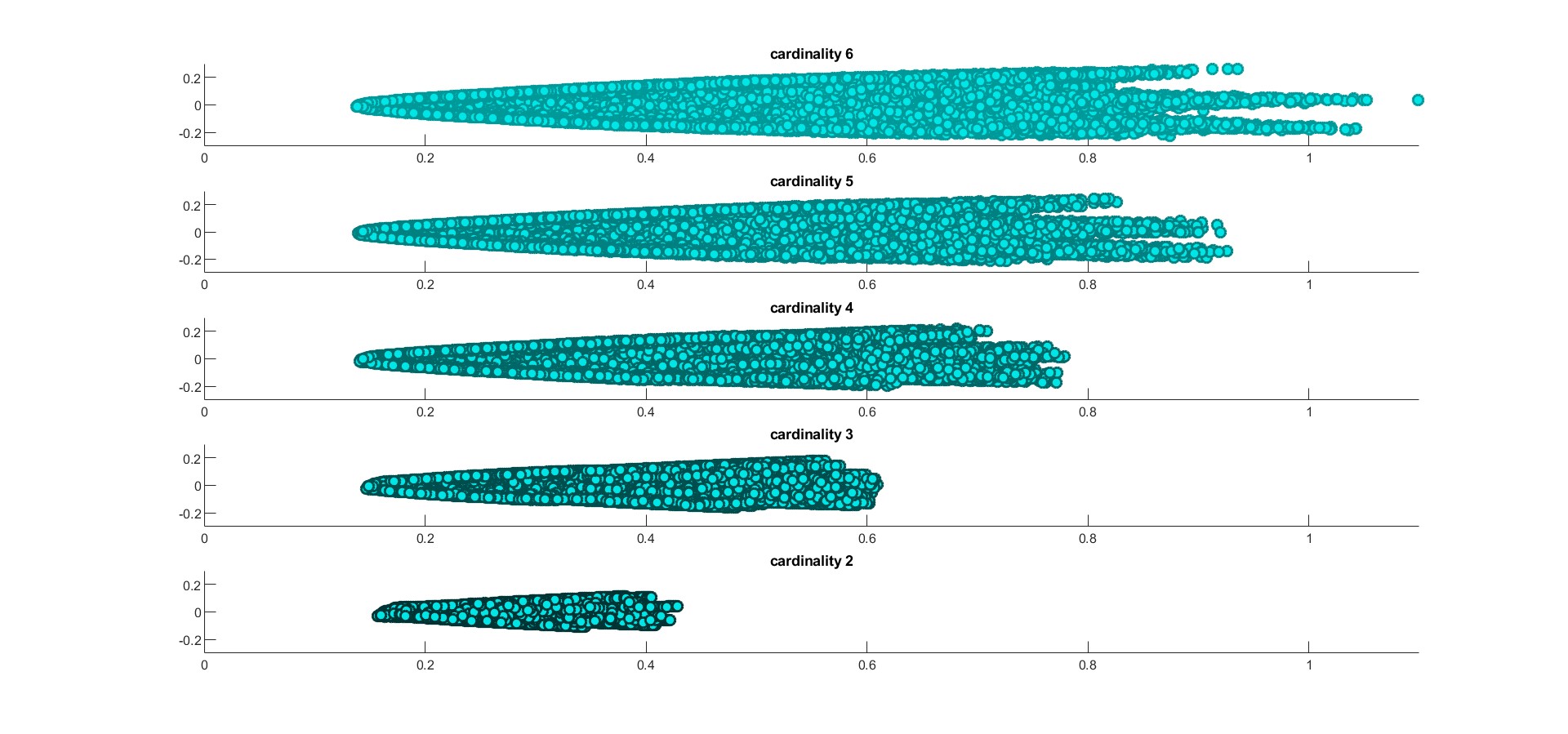}
	\end{center}
	\caption{Feasible set for the Simple case and  $\alpha=2,3,4,5$ and $6$.}
	\label{fig18}
\end{figure}

 For the Simple-case problem, with $n=6$ available assets, the feasible set
  is shown in Figure \ref{fig18}. We note that for larger value of  $\alpha$ we get  larger  area of  the feasible set.
 We also note that the bullet shape is not affected by the cardinality constraint but, as expected, the set is shrinking as the number of
  zero elements increases.

\begin{figure}[H]
	\begin{center}\hspace{-2cm}
		\includegraphics[width=19cm]{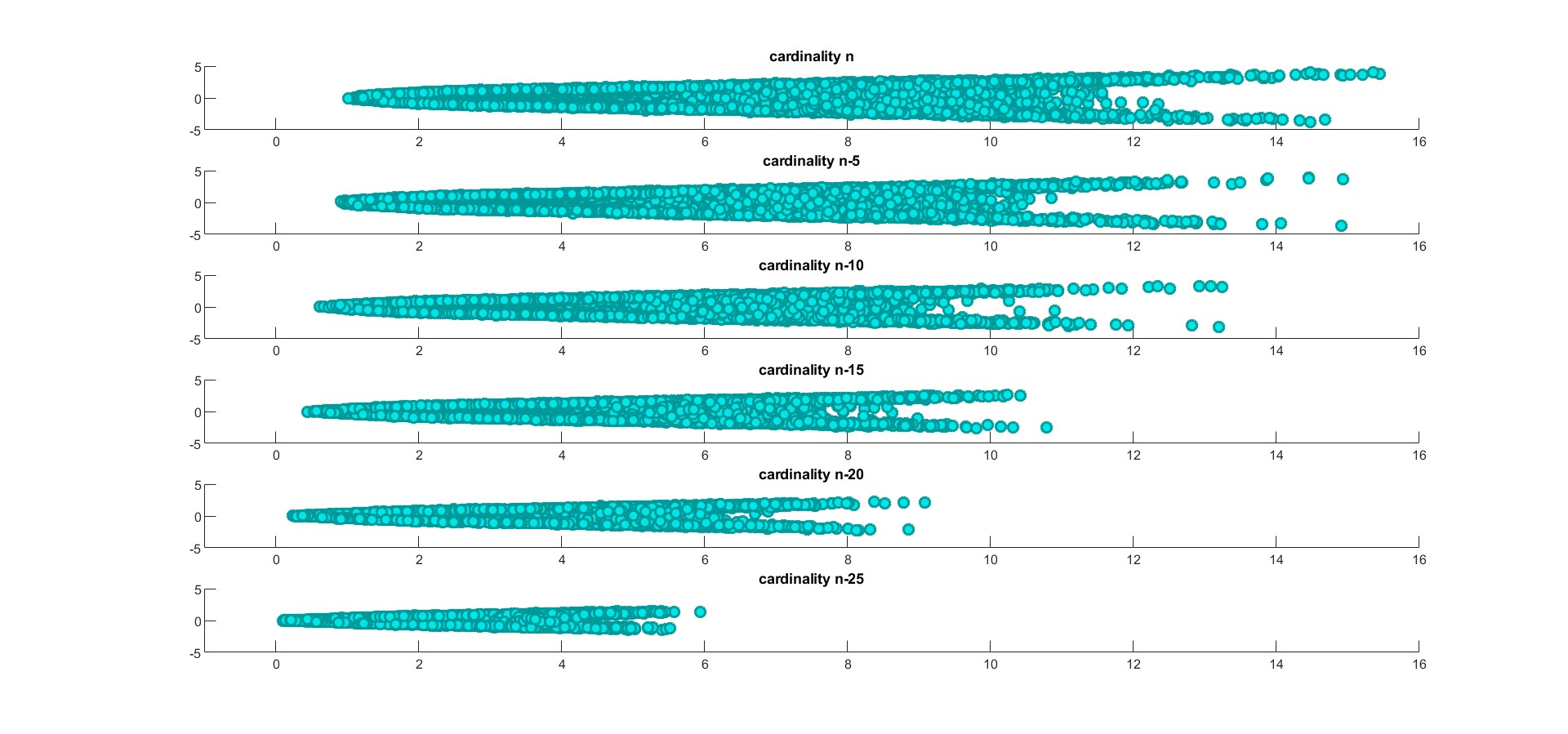}
	\end{center}
	\caption{Feasible set for  Port1 and  $\alpha=6, 11, 16, 21, 26$ and $31$.}
	\label{fig21}
\end{figure}

\begin{figure}[H]
	\begin{center}\hspace{-2cm}
		\includegraphics[width=19cm]{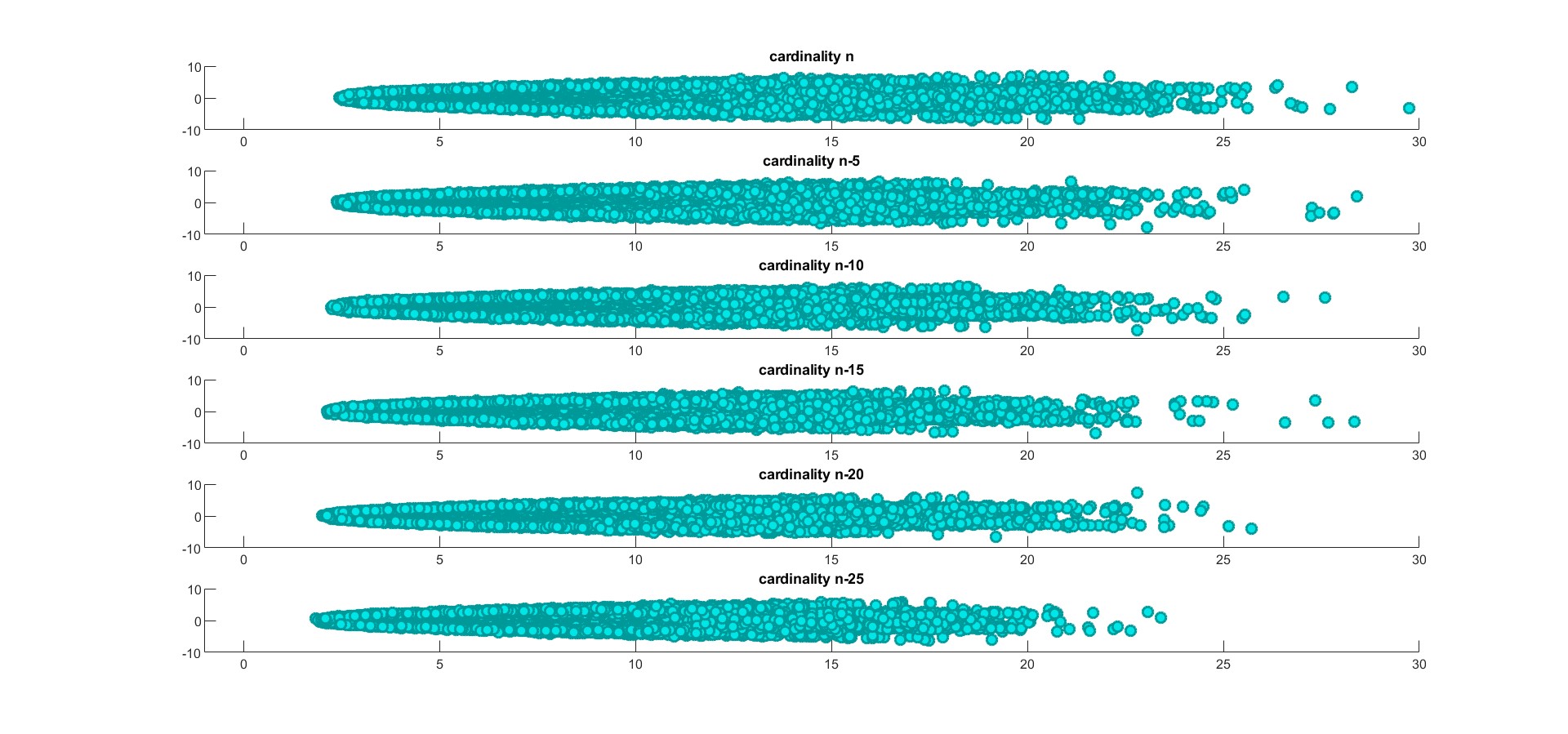}
	\end{center}
	\caption{Feasible set for  Port2 and  $\alpha=60, 65, 70, 75, 80$ and $85$.}
	\label{fig22}
\end{figure}

The same conclusions apply to the larger data sets coming from real assets. Below, in Figures \ref{fig21} and \ref{fig22}, we show the feasible sets
for Port1 and Port 2. We note that once again the area is shrinking when $\alpha$ decreases. We also note that the same is true for all considered cases.

\begin{figure}[H]
	\begin{center}
		\includegraphics[width=18cm]{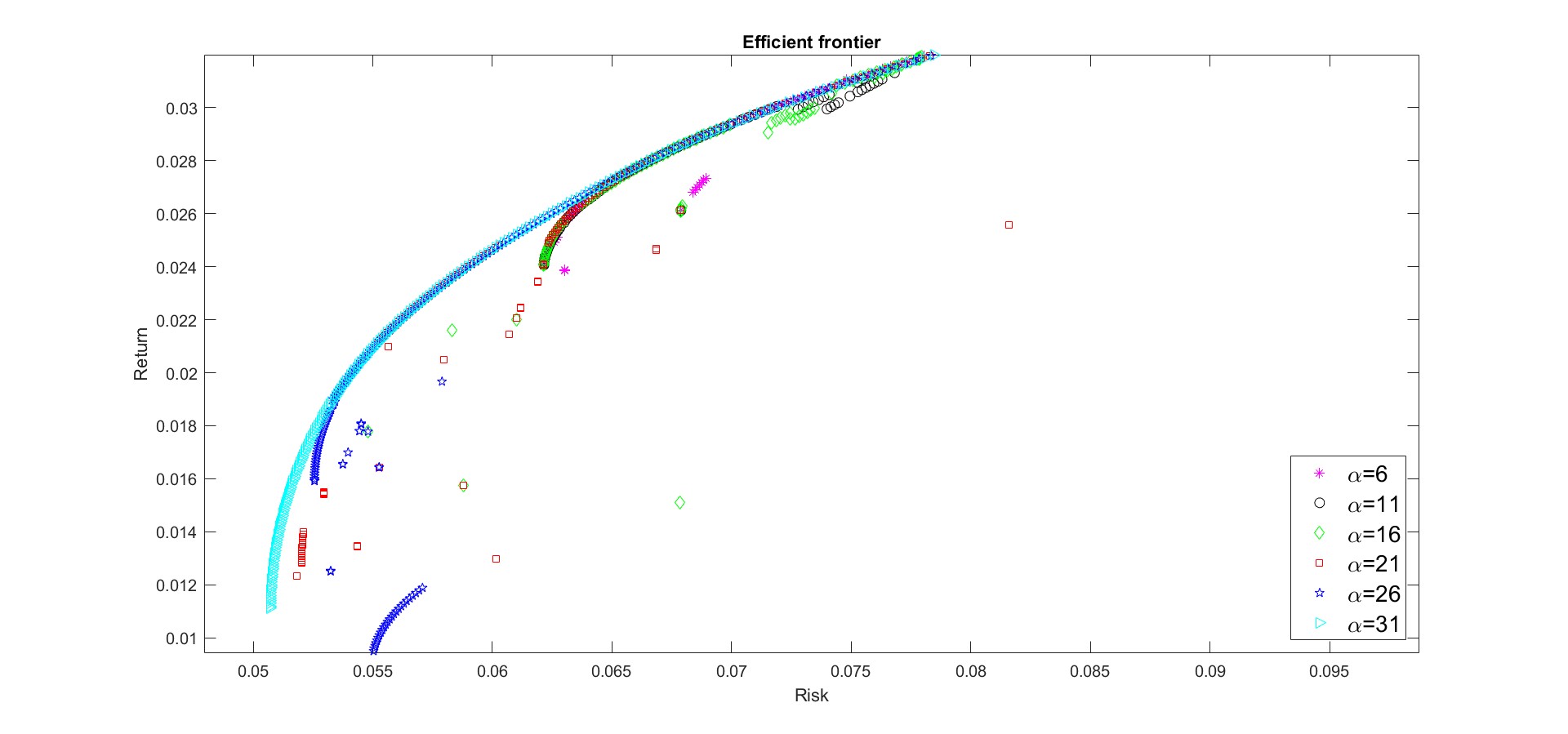}
	\end{center}
	\caption{Risk versus return, using Algorithm PSPGD for Port1 and $\alpha=6,11,16,21,26, 31$.}
\label{fig5}
\end{figure}

\begin{figure}[H]
	\begin{center}
		\includegraphics[width=18cm]{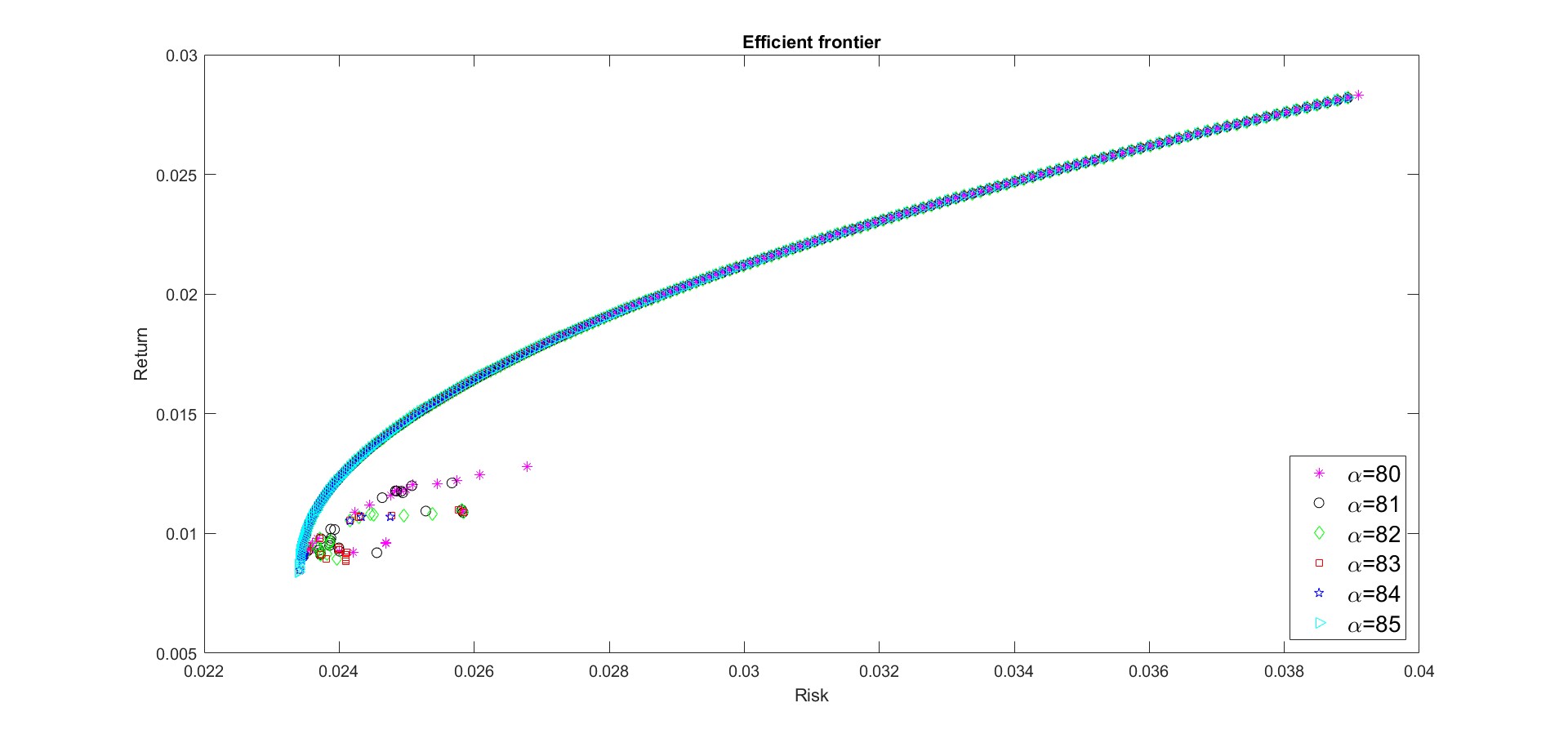}
	\end{center}
	\caption{Risk versus return, using Algorithm PSPGD for  Port2 and $\alpha=80,81,82,83,84,85$. }
\label{fig10}
\end{figure}

\begin{figure}[H]
	\begin{center}
		\includegraphics[width=18cm]{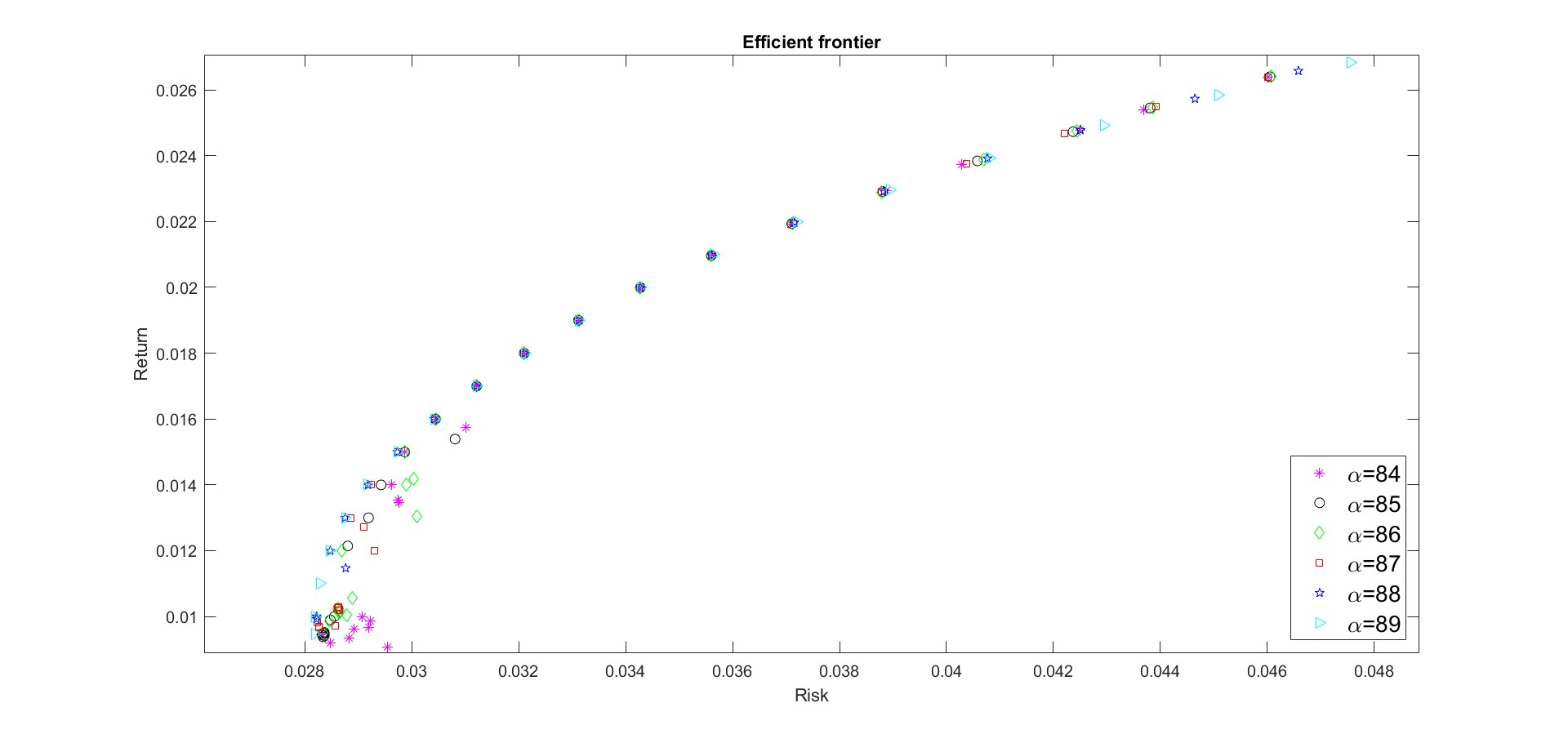}
	\end{center}
	\caption{Risk versus return, using Algorithm PSPGD for Port3 and $\alpha=84,85,86,87,88,89$. }
\label{fig27}
\end{figure}

\begin{figure}[H]
	\begin{center}
		\includegraphics[width=18cm]{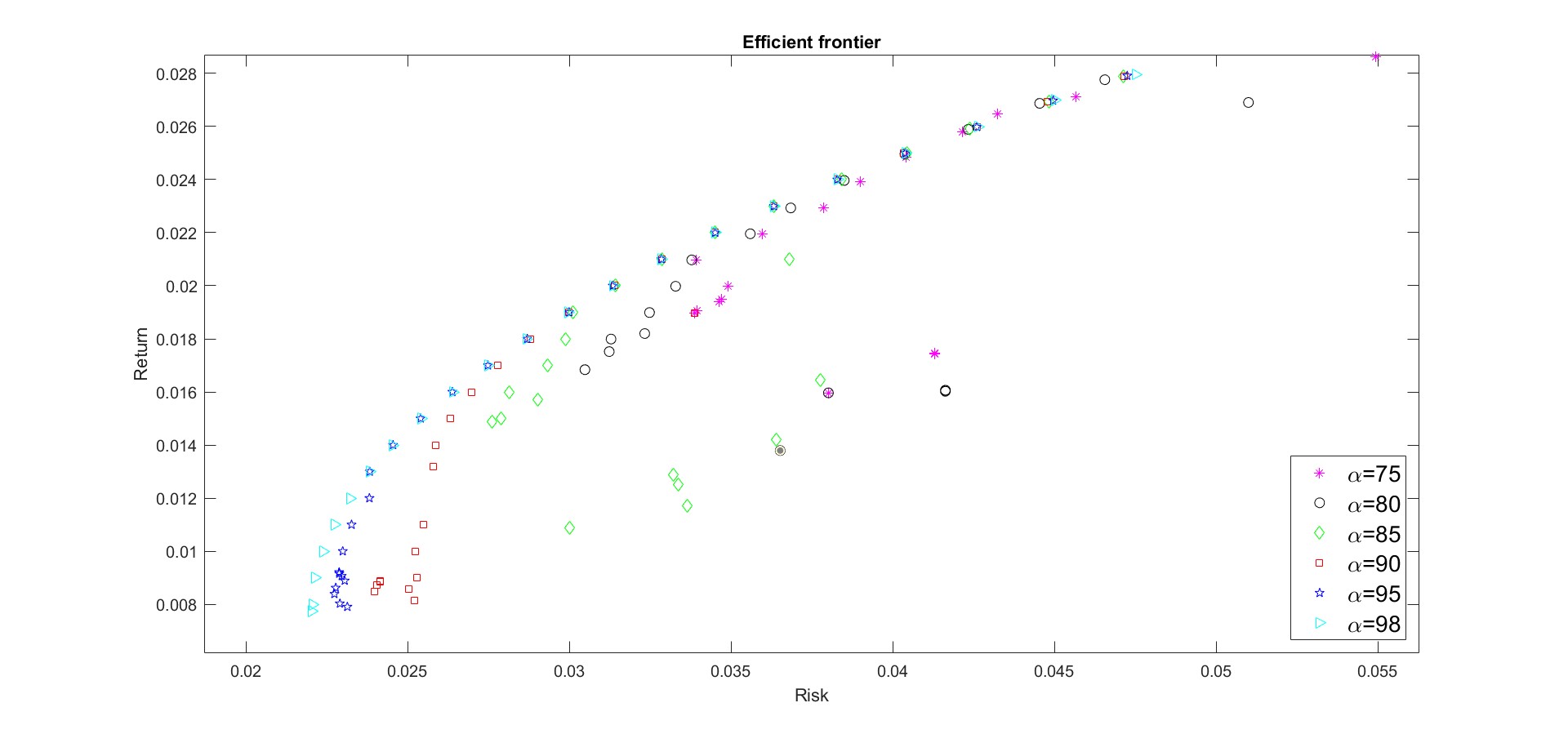}
	\end{center}
	\caption{Risk versus return, using Algorithm PSPGD for  Port4 and  $\alpha=75,80,85,80,95,98$. }
\label{fig28}
\end{figure}

\begin{figure}[H]
	\begin{center}
		\includegraphics[width=18cm]{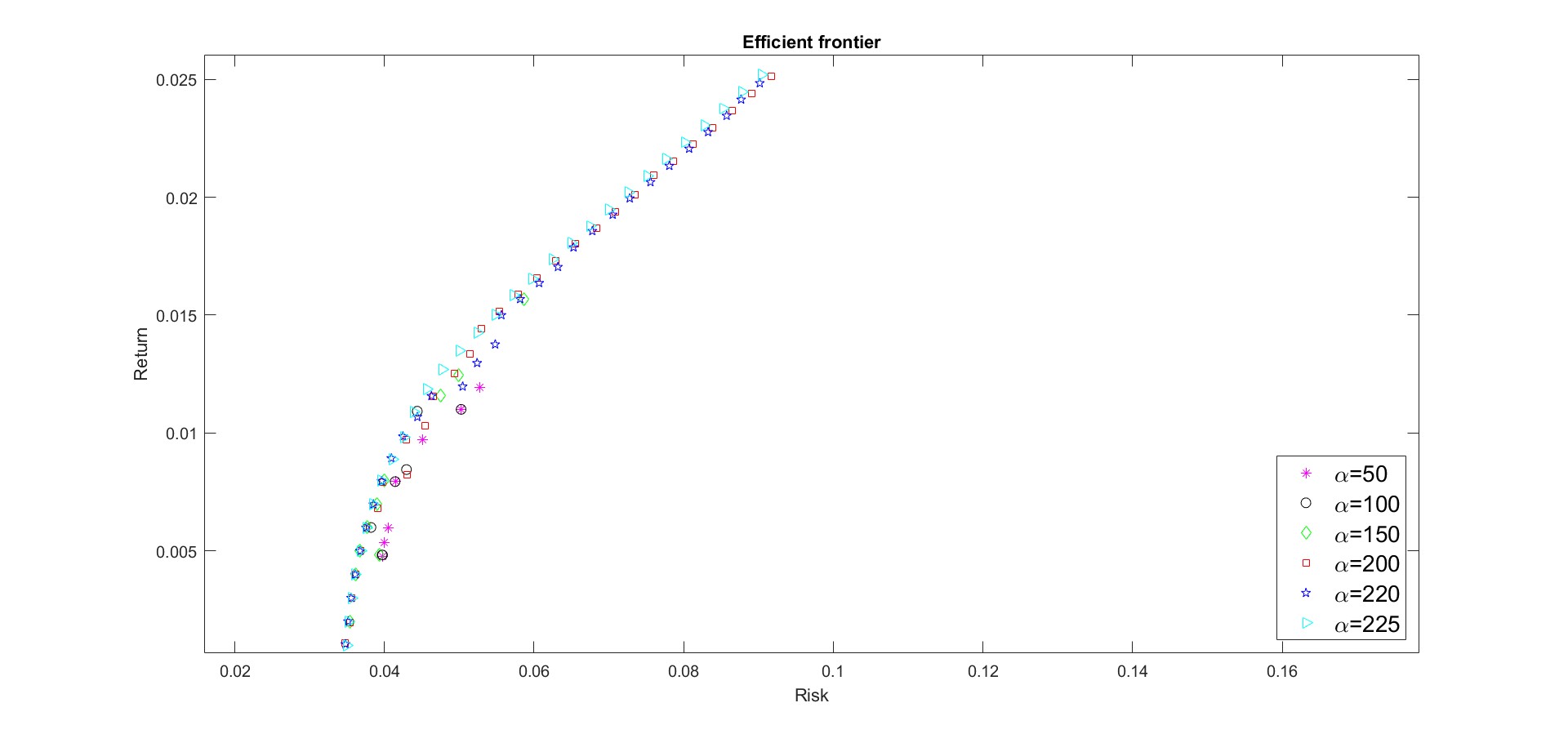}
	\end{center}
	\caption{Risk versus return, using Algorithm PSPGD for  Port5 and $\alpha=50,100,150,200,220,225$.}
\label{fig29}
\end{figure}

\begin{figure}[H]
	\begin{center}
		\includegraphics[width=18cm]{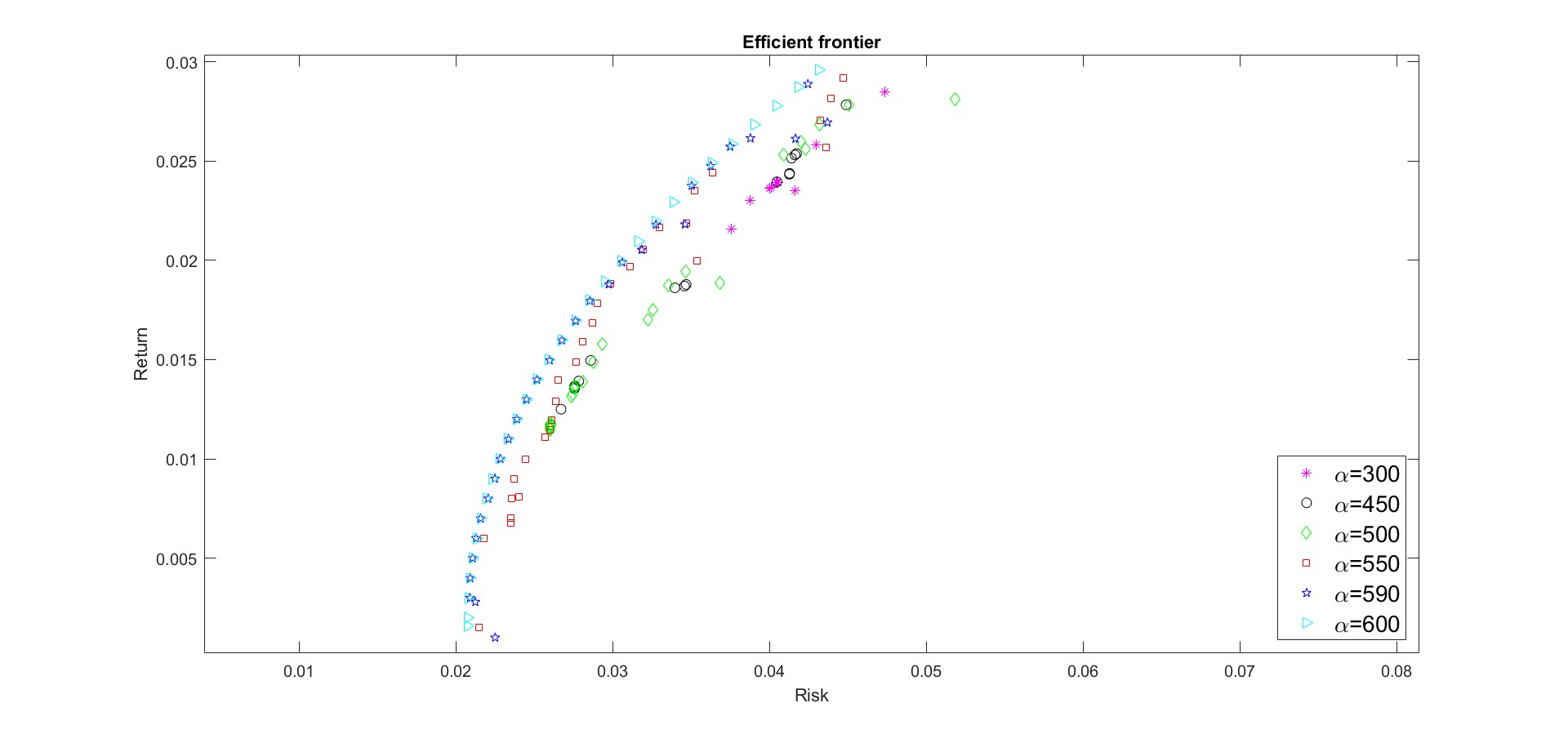}
	\end{center}
	\caption{Risk versus return, using Algorithm PSPGD for Port6 and $\alpha=300,450,500,550,590,600$.}
\label{fig30}
\end{figure}

The efficient frontiers for all data sets are shown in Figures \ref{fig5}--\ref{fig30}. Again, we observe that the efficient frontier
 is deformed by the value of the cardinality constraint, and when $\alpha < n$ it is not a convex curve.  For the sake of
 completeness, in the Appendix we provide some tables with more detailed results, varying the cardinality constraints, for all considered data sets.
 We can observe in all those figures and tables the effectiveness of our low-cost continuous approach (Algorithm PSPGD).

Additionally we compare our approach to IBM ILOG CPLEX Optimization Studio,  Version: 22.1.0.0. CPLEX is a  quadratic mixed integer programming solver.
 The goal of comparison is to investigate the quality of solutions obtained by PSPGD and CPLEX in term of risk and return. We also report CPU time 
 although CPLEX is implemented in a low-level language and so it requires significantly less execution time than our high-level Matlab implementation. 
 Hence, CPU time might be misleading.  For solving the problems with CPLEX we consider the following formulation:
	\begin{equation*}
		\begin{array}{ll}
			\displaystyle	\min_{x,y}   & \frac{1}{2} x^{\top}Qx\\ [1mm]
			\mbox{subject to:} & v^{\top} x \geq \rho,  \\ [1mm]
			& e^{\top} x = 1,  \\ [1mm]
			& e^{\top} y \geq n-\alpha,  \\ [1mm]
			& 0 \leq x_i\leq 1,  \; \mbox{ for all } 1\leq i \leq n,  \\ [1mm]
			& x_i+y_i\leq 1, 	 \; \mbox{ for all } 1\leq i \leq n,  \\ [1mm]
			&  y_i \in \{0,1\}.
		\end{array}
	\end{equation*}
Notice that in the above problem formulation we do not have the Hadamard constrained and instead  we have $ x_i + y_i \leq 1 $ followed by $ y_i \in  \{0,1\}$. 
 CPLEX is designed to work with linear constraints and for $ y_i = 0 $ or $ y_i = 1 $ we get the same condition.

The details of tests for all considered data sets are presented in Tables \ref{Table:simple} - \ref{Table:port6} in Appendix. Notice that CPLEX data is missing for Port 6 as we were not able to solve the problem with CPLEX.
One can easily see that PSPGD produces solutions with slightly higher risk  and significantly better return. 	
 In Tables~\ref{Table:port2_b}  we observe that CPLEX needs a very large number of iterations to solve the problem for  $\alpha\leq 20$, which corresponds to the fact the PSPGD needed a special value of  $\tau_{-1}$ for these values of  $\alpha$ and large values of penalty parameter $ \tau. $ Thus, this behavior is associated with the data of Port2.  In some other cases, reported in the tables in Appendix, we can observe a rather large number of CPLEX iterations for
  small values of $\alpha$ while PSPGD solved the same problems with  reasonably small values of the penalty parameters.

An interesting observation from the literature, and confirmed by our experiments, is the fact that the optimal portfolio without cardinality
constraint is in fact sparse. In Table~\ref{Tabl8} we report the number of assets obtained by our algorithm and CPLEX which is in accordance with the results
reported in  \cite[Figure 5]{Cesarone2009} and \cite[Section 5.2.2]{Cesarone2013}. We can observe that  the number of assets
in the unconstrained Mean-Variance optimal portfolio for Port1 $\|x^*\|_0\leq 12$, for Port4 $\|x^*\|_0\leq 40$,  and for Port5 $\|x^*\|_0\leq 15$.

\begin{table}[H]
	\caption{Performance of Algorithm PSPGD  for all cases when $n=\alpha$.}  \label{Tabl8}
	\begin{center}
		\begin{tabular}{|c|c|ccc|ccc|}\hline
			&&\multicolumn{3}{c|}{PSPGD}&\multicolumn{3}{c|}{CPLEX}\\
			Problem&	$\alpha=n$ &  $\|x\|_0$& $v^{\top} x^*$ & $\sqrt{(x^*)^{\top}Qx^*}$   &  $\|x\|_0$& $v^{\top} x^*$ & $\sqrt{(x^*)^{\top}Qx^*}$    \\
			\hline \hline
			Simple-case &	6   & 6 & 0.0003& 0.1394 &6&0.0003&0.1394\\
			Port1 &	31   & {12} & {0.0133}& {0.0509}&12&0.0133&0.0509 \\
			Port2 &	85   & 	{24} & {0.0085}& 0.0234 &25
			&0.0084&0.0234\\		
			Port3 &	89   & {34} & 0.0101& 0.0282 &34&0.0101&0.0282\\		
			Port4 &	98   & {38} & {0.0098}& {0.0223} &38&0.0098&0.0223\\		
			Port5 &	225   & {12} &{0.0003} &{0.0349} &12 &0.0003&0.0349\\	
		    Port6 &	600   & 38 & 0.0015& 0.0207 &&&\\	\hline				\end{tabular}
\end{center}
\end{table}

	As noticed above the feasible set of (\ref{portfol}) belongs to the feasible set of (\ref{portfol2}). In addition, since the solution of
(\ref{portfol2}) satisfies the Hadamard condition we obtain that the solution is also a solution of (\ref{portfol}). Then, by
Theorem~\ref{theo: min_bin}, we have that if $({x}^*,{y}^*)$ is a local minimizer of (\ref{portfol2}) satisfying $\|{x}^*\|_0=\alpha$ then the components of ${y}^*$ are binary, ${y}^*$  is unique, and  $x^*$ is a local minimizer of (\ref{origprb}). In fact, for the solutions reported in
 Tables~\ref{Table:simple},~\ref{Table:port1}, and~\ref{Table:port3_b}  in the Appendix, if $\|x^*\|_0=\alpha$ we have that
  the components of ${y}^*$ are binary. The solution may have non-binary entries in $y^*$, for instance port1 with $\alpha=n=31$ we have that $y^*$
    is binary, however the cardinality constraint is not active $\|x^*\|_0=12$. Another interesting example is detected for Port3 with $\alpha=n=89$ in which we obtain
    a binary $y^*$ but $\|x^*\|_0=34$.

\section{Conclusions and final remarks} \label{conclu}

Taking advantage of a recently developed continuous formulation, we have developed and analyzed a low-cost and effective scheme for
solving convex constrained optimization problems that also include a ``hard-to-deal" cardinality constraint.
As it appears in many applications, we assume that the region defined by the convex constraints
 can be written as the intersection of a finite collection of ``easy to project" convex sets.
Under this continuous formulation, to fulfill the cardinality constraint, the Hadamard condition $x\circ y  = 0$
 must be satisfied between the solution vector $x$ and an auxiliary vector $y$. In our scheme this
  condition is achieved by adding a non-negative penalty term $h(x,y)$,  and using a classical penalization
  strategy. For each penalty subproblem, a convex constrained
  problem must be solved, which in our proposal is achieved by combining two low-cost computational schemes:
  the spectral projected gradient (SPG) method  and   Dykstra's alternating projection method. \\

  To illustrate the computational performance of our combined scheme, we have considered in detail the standard
    mean-variance limited diversified portfolio selection problem, which involves obtaining the proportion of the initial budget that should be
    allocated in a limited number of the available assets. For this specific application, we proposed a natural differentiable choice of the penalty
     term (given by $h(x,y)= x^{\top}y$) that must be driven to zero, which  allowed us to develop a simple way of increasing the associated penalty
      parameter in a controlled and bounded way.  In our numerical study we have included a variety of data sets involving real-world capital
       market indices. For these data sets we have produced the feasible sets and also the efficient frontier (a curve illustrating the tradeoff
 between risk and return) for different values of the limited number of allowed assets. In each case, we  highlighted the differences that arise
 in the shape of this efficient frontiers as compared  with the unconstrained efficient one. The presented numerical study includes comparison with CPLEX, 
 a professional software for general mixed integer programming problems.  The comparison is presented in terms of quality of solution (higher return, 
 lower risk) and PSPGD appears to be competitive. Furthermore, PSPGD is successfully applied to a large portfolio problem with 600 assets while CPLEX 
 failed at solving this particular problem. \\

  In our modeling of the portfolio problem we have bounded the proportion to be invested in each of the selected assets between 0 and 1. However,
   without altering our proposed scheme, stricter upper limits (less than 1) can
   be imposed on some particular assets. Clearly, this would require a more careful analysis of the feasible
   options for the expected return. Moreover, it could also be interesting from a portfolio point of view,
   to allow negative entries in some  of the proportions to be invested, and that can be accomplished by
   allowing negative values in the  lower bounds of the solution vector. In that case, the penalization
   term to force the Hadamard condition needs to be chosen accordingly (e.g., $h(x,y) = \sum_{i=1}^n (x_i^2 y_i)$). \\ [2mm]

\noindent
 {\bf Acknowledgements.} \\
   The first author was financially supported by the  Serbian  Ministry  of Education, Science, and Technological Development and Serbian Academy of Science and Arts,
  grant no.  F10. The second author was financially supported by  Funda\c c\~ao para a Ci\^encia e a Tecnologia (FCT) (Portuguese Foundation for Science
  and Technology) under the scope of the projects 	UIDB/MAT/00297/2020, UIDP/MAT/00297/2020
  (Centro de Matem\'atica e Aplica\c{c}\~oes), and UI/297/2020-5/2021. The third author  was financially supported by  the Funda\c c\~ao para a Ci\^encia e a Tecnologia 
  (Portuguese Foundation for Science and Technology) 
  under the scope of the projects UIDB/MAT/00297/2020 and  UIDP/MAT/00297/2020 (Centro de Matem\'atica e Aplica\c{c}\~oes). \\ [2mm]

\noindent
{\bf Data availability.}  The codes and data sets generated during and/or analyzed during the current
study are available from the corresponding author on reasonable request. \\ [2mm]

\noindent
{\bf Disclosure statement.} No potential conflict of interest was reported by the authors. \\ [2mm]

\section*{Appendix: Performance of Algorithm PSPGD for all data sets}

In Tables~\ref{Table:simple} - \ref{Table:port5_b}, we report the performance of PSPGD and CPLEX, for several values of $\alpha$,
 reporting the values of optimal portfolio return, risk, number of non-zero portfolio weights,  number of iteration (Iter), and number of SPG
  iterations for PSPGD,  the CPU time (Time) in seconds,  the last value of $\tau$, as well as the final value of the Hadamard product,
 and the total number of required function evaluations (fcnt).  It is worth noticing that in all the results reported in these tables,
the  pgnorm at the obtained solution and the Hadamard products  $(x^*)^{\top}y^*$ are strictly less than $10^{-6}$, and hence we did not report these values.  In Table \ref{Table:port6} the results rae reported only for PSPGD as we were not able to solve Port 6 with CPLEX.

\begin{table}[H]
	\caption{Performance of PSPGD and CPLEX for the Simple case.}  \label{Table:simple}
	\begin{center}
		\begin{tabular}{|ccccccccccc|}\hline
			Algorithm  &$\alpha$ & $v^{\top}x^*$ &$\sqrt{(x^*)^{\top}Qx^*}$ & $\|x^*\|_0$ & Iter & Iter-SPG & Time & $\tau$  & fcnt & $\rho$\\ \hline \hline PSPGD
&  1  & 0.0400 &0.2074&1&   2 & 4 & 0.3708 &  0.117590 & 7 & 0.0018\\
& 2 &  0.0293 &0.1735&2 & 2 & 6 &0.3133 & 0.117577 & 9  &  0.0016 \\
& 3 & 0.0053 &0.1523&3&   2 & 11 &  0.2786 &0.117560 & 13 & 0.0017  \\
 & 4  & 0.0053 &0.1523&3& 2 & 12 & 0.3211 & 0.117558 &16 & 0.0017 \\
 &  5  & 0.0053 &0.1523&3&  2 &  8 & 0.2799 &  0.117557  &10 & 0.0012 \\
 &6  & 0.0003 &0.1394&6&  2 & 7 &  0.3001  &  0.117556 &9  & 0.0003\\ \hline \hline
 CPLEX & 1 &  0.0210&0.1949&1&22 & - &0.09 & - & -  & 0.0018\\
 & 2 &0.0016&0.1612&2&19 &- &  0.05 & - & - & 0.0016 \\
 & 3 & 0.0017&0.1483&3&19 & - & 0.03&  - & - & 0.0017\\
 & 4 & 0.0017 & 0.1414 & 4 & 19 & - & 0.05 &  - & - & 0.0017 \\
 & 5 & 0.0012 & 0.1414 & 5 & 19 & - & 0.06 & - & - & 0.0012 \\
 & 6 & 0.0003 & 0.1394 & 6 & 13 & - & 0.02 & - & - & 0.0003 \\
 \hline
\end{tabular}
\end{center}
\end{table}

\begin{table}[H]
	\caption{Performance of Algorithm PSPGD and CPLEX for problem Port1.}  \label{Table:port1}
	\begin{center}
		\begin{tabular}{|ccccccccccc|}\hline
			Algorithm  &$\alpha$ & $v^{\top}x^*$ &$\sqrt{(x^*)^{\top}Qx^*}$ & $\|x^*\|_0$ & Iter & Iter-SPG & Time & $\tau$  & fcnt & $\rho$\\ \hline \hline PSPGD
&1 & 0.0435  &  0.1382 &         1    & 2 &  4 & 0.8109  &  0.1475 &    6& 0.0097\\
&2 & 0.0435  &  0.1382 &         1    & 2 & 3 & 0.3113  &  0.1476 &      5&0.0126 \\
&3 & 0.0435  &  0.1382 &         1    & 2 & 3 & 0.3203  &  0.1477 &      5 &0.0133\\
&4 & 0.0435  &  0.1382 &         1    & 2 & 3 & 0.2200  &  0.1476 &       5 &0.0132\\
&5 & 0.0435  &  0.1382 &         1    & 2 & 3 & 0.2524  &  0.1476 &      5 &0.0133\\
&10 & 0.0435 &  0.1382 &         1    &  2 & 4 & 0.2028  &  0.1475 &      6 &0.0136\\
&15 & 0.0151 &  0.0678 &         2    & 2 & 17 & 0.6132  &  0.1473 &      23  & 0.0133\\
&20 & 0.0154 &  0.0530 &         5    & 2 & 17 & 0.3751  &  0.1473 &      19  & 0.0132\\
&30 & 0.0133 &  0.0509 &        11    & 2 & 13 & 0.2978  &  0.1471 &      15  & 0.0133\\
&31 & 0.0133 &  0.0509 &         12 &  2 & 12 & 0.3267  &  0.1471 &       14  & 0.0133\\\hline \hline
CPLEX & 1 & 0.0233 & 0.0717&1& 32& --&0.0900&--&--&0.0097\\
 & 2 & 0.0126  & 0.0591&2& 17& --&0.0300&--&--&0.0126\\
 & 3 & 0.0140  & 0.0544&3& 17& --&0.0500&--&--&0.0133\\
 & 4 & 0.0132 & 0.0523&4& 17& --&0.0300&--&--&0.0132\\
 & 5 & 0.0137 & 0.0516&1&17&--&0.0500&--&--&0.0133\\
 & 10 & 0.0136& 0.0510&10& 19& --&0.0600&--&--&0.0136\\
  & 15 & 0.0133 & 0.0509&12&13&--&0.0300&--&--&0.0133\\
 & 20 & 0.0132 & 0.0509&12& 13& --&0.0200&--&--&0.0132\\
 & 30 & 0.0133 & 0.0509&12&13&--&0.0200&--&--&0.0133\\
 & 31 & 0.0133 & 0.0509&12&13&--&0.0200&--&--&0.0133\\   \hline
		\end{tabular}
\end{center}
\end{table}

\begin{table}[H]
	\caption{Performance of Algorithm  PSPGD and CPLEX for problem Port2.}  \label{Table:port2_b}
	\begin{center}
		\begin{tabular}{|ccccccccccc|}\hline
			Problem&$\alpha$ & $v^{\top}x^*$ &$\sqrt{(x^*)^{\top}Qx^*}$ & $\|x^*\|_0$ & Iter & Iter-SPG & Time & $\tau$  & fcnt & $\rho$\\ \hline \hline
			PSPGD &  1&0.0392  &0.1065&1&2&8&0.6761&0.0976&28&0.0085 \\
 &      2 & 0.0392 &  0.1065 &         1    &          2  &   5 & 0.3885  &  7.0147 &  19 & 0.0058 \\
&       3 & 0.0745 &  0.1327 &         2    &          2  &  10 & 0.5509  &  11.815 &  39 & 0.0079 \\
&4 & 0.1045 &  0.1628 &         3    &          3  &  29 & 0.7973  &  12.012 &   155 & 0.0125 \\
&5 & 0.0745 &  0.1327 &         2    &          2  &  10 & 0.5157  &  11.866 &  42 & 0.0163
\\
     & 10 & 0.1267 &  0.2010 &         4    &          5  & 109 & 2.5451  &  42.404 &  383 & 0.0158 \\
       &    15 & 0.1804 &  0.2954 &         7    &          3  &  58 & 1.4936  &  72.623 &    163 & 0.0161 \\
     &  20 & 0.0745 &  0.1327 &         2    &          2  &  12 & 0.5987  &  26.025 &   69 & 0.0022 \\
  & 25 & 0.0745 &  0.1327 &         2    &          2  &  12 & 2.1676  &  29.028 &    69 &0.0024\\
        &30 & 0.0745 &  0.1327 &         2    &          2  &  11 & 0.8536  &  39.052 &   61 & 0.0037\\
        &	35	&	0.0291	&	0.0428	&	5	&	2	&	11	&	0.2984	&	0.0977	&	18	&	0.0109	\\
&	40	&	0.0291	&	0.0428	&	5	&	2	&	11	&	0.2845	&	0.0977	&	18	&	0.0117	\\
&	45	&	0.0291	&	0.0428	&	5	&	2	&	13	&	0.2991	&	0.0976	&	25	&	0.0115	\\
&	50	&	0.0291	&	0.0428	&	5	&	2	&	11	&	0.2731	&	0.0976	&	18	&	0.0115	\\
&	55	&	0.0225	&	0.0357	&	8	&	2	&	15	&	0.2739	&	0.0975	&	22	&	0.0110	\\
&	60	&	0.0186	&	0.0319	&	13	&	2	&	13	&	0.2580	&	0.0974	&	17	&	0.0110	\\
&	65	&	0.0190	&	0.0321	&	12	&	2	&	15	&	0.2748	&	0.0974	&	17	&	0.0111	\\
			&	70	&	0.0110	&	0.0237	&	23	&	2	&	19	&	0.2579	&	0.0973	&	21	&	0.0110	\\
&	75	&	0.0111	&	0.0238	&	23	&	2	&	16	&	0.2120	&	0.0973	&	18	&	0.0111	\\
&	80	&	0.0103	&	0.0235	&	25	&	2	&	16	&	0.2108	&	0.0973	&	18	&	0.0103	\\
&	85	&	0.0085	&	0.0234	&	24	&	2	&	15	&	0.1876	&	0.0973	&	17	&	0.0070	\\
\hline \hline
CPLEX
& 1 &0.0134&0.0477 &1&504 & --&0.14&--&--&0.0085\\
& 2 &0.0066 &0.0331 &2&5638 & --&0.27&--&--&0.0058\\
& 3 & 0.0084 &0.0296 &3& 34024& --&0.53&--&--&0.0079\\
& 4 &  0.0125& 0.0289&4&13926 & --&0.39&--&--&0.0125\\
& 5 &0.0163  &0.0298 &5& 6982& --&0.28&--&--&0.0163\\
& 10 & 0.0158 &0.0263 &10& 2743& --&0.20&--&--& 0.0158\\
& 15 & 0.0161 & 0.0259&15& 1485& --&0.23&--&--&0.0161\\
& 20 &0.0083  &0.0234 &20& 75& --&0.20&--&--&0.0022\\
& 25 & 0.0083 &0.0234 &25& 13& --&0.02&--&--&0.0024\\
& 30 &0.0084  &0.0234 &25& 13& --&0.02&--&--&0.0037\\
& 35 &0.0109  &0.0236 &24&14 & --&0.02&--&--&0.0109\\
& 40 & 0.0117 & 0.0238&24&14 & --&0.03&--&--&0.0117\\
& 45 & 0.0115 &0.0238 &24&14 & --&0.03&--&--&0.0115\\
& 50 & 0.0115 &0.0238 &24& 14& --&0.02&--&--&0.0115\\
& 55 &0.0110  &0.0237 &24&14 & --&0.02&--&--&0.0110\\
& 60 & 0.0110 &0.0237 &24&14 & --&0.02&--&--&0.0110\\
& 65 &0.0111  & 0.0237&24&14 & --&0.02&--&--&0.0111\\
& 70 & 0.0110 &0.0237 &24&14 & --&0.03&--&--&0.0110\\
& 75 & 0.0111 & 0.0237&24&14 & --&0.05&--&--&0.0111\\
& 80 & 0.0103  & 0.0235&26&14 & --&0.03&--&--&0.0103\\
& 85 & 0.0084 & 0.0234& 25 & 13& --&0.03&--&--&0.0070\\
\hline
		\end{tabular}
	\end{center}
\end{table}

\begin{table}[H]
	\caption{Performance of Algorithm PSPGD and CPLEX for problem Port3.}  \label{Table:port3_b}
	\begin{center}
		\begin{tabular}{|ccccccccccc|}\hline
			Problem&$\alpha$ & $v^{\top}x^*$ &$\sqrt{(x^*)^{\top}Qx^*}$ & $\|x^*\|_0$ & Iter & Iter-SPG & Time & $\tau$  & fcnt & $\rho$\\ \hline \hline
			PSPGD&       1 & 0.0328 &  0.0779 &         1    &          2  &   7 & 1.0909  &  0.1133  &   16 & 0.0101 \\		 	
		& 2 & 0.0328 &  0.0779 &         1    &          2  &   3 & 1.0322  &  9.8615 &    11 & 0.0104 \\
			&3 & 0.0328 &  0.0779 &         1    &          2  &   4 & 0.6540  &  9.7468 &     12 & 0.0102 \\
			&4 & 0.0328 &  0.0779 &         1    &          2  &   3 & 0.5711  &  9.6557 &     10 & 0.0160 \\
			&5 & 0.0328 &  0.0779 &         1    &          2  &   5 & 0.7342  &  9.5319 &     17 & 0.0135 \\
&10 & 0.0328 &  0.0779 &         1    &          3  &  56 & 4.5180  &  9.0391 &  402 & 0.0119 \\
&15 & 0.0328 &  0.0779 &         1    &          3  &  55 & 1.3683  &  11.429  &  406 & 0.0117 \\
  &   20 & 0.0104 &  0.0284 &        14    &         14  & 676 & 36.438  &  0.0003  &  885 & 0.0104 \\
&25 & 0.0104 &  0.0284 &        14    &         11  & 530 & 28.613  &  0.0003  &  666 & 0.0104 \\
&30 & 0.0104 &  0.0284 &        14    &         12  & 547 & 29.572  &  0.0003  &  695 & 0.0104 \\
 &35 & 0.0328 &  0.0779 &         1    &          2  &   4 & 0.5872  &  89.248  &   13 & 0.0107 \\
      &40 & 0.0104 &  0.0286 &        12    &          5  & 215 & 10.943  &  0.0005 &  282 & 0.0104 \\
&       45 & 0.0251 &  0.0464 &         4    &          3  &  51 & 1.9333  &  0.2019 &   59 & 0.0114 \\
&50 & 0.0104 &  0.0285 &        13    &          5  & 194 & 3.2324  &  0.0005 &   252 & 0.0104 \\
&55 & 0.0104 &  0.0285 &        13    &          4  & 151 & 2.6642  &  0.0005 &  183 & 0.0104 \\
&60 & 0.0167 &  0.0321 &        19    &          5  & 183 & 6.3165  &  1.1546  &  440 & 0.0133 \\
&65	& 0.0157 	&  0.0310& 	21&	8&	341&	 5.8265&  	  2.8615 	&652	& 0.0145 \\
&	70	&	 0.0105 	&	  0.0333 	&	10	&	2	&	29	&	 0.9089  	&	  0.1129 	&	31	&	 0.0105 	\\
&	75	&	 0.0105 	&	  0.0308 	&	11	&	2	&	20	&	 0.6148  	&	  0.1129 	&	22	&	 0.0105 	\\
&	80	&	 0.0105 	&	  0.0295 	&	16	&	2	&	16	&	 0.6074  	&	  0.1129 	&	18	&	 0.0105 	\\
&	85	&	 0.0104 	&	  0.0286 	&	23	&	2	&	26	&	 0.6432  	&	  0.1129 	&	28	&	 0.0104 	\\
&	89	&	 0.0101 	&	  0.0282 	& 34	&	2	&	18	&	 0.5024  	&	  0.1129 	&	20	&	 0.0101	\\
\hline   \hline
CPLEX
& 1 & 0.0151&0.0473 &1 &328 & --&0.19 &--&--&0.0101\\
& 2 &0.0117 &0.0384 &2 &9537 & --&0.42&--&--&0.0104\\
& 3 &0.0104  & 0.0346&3&133879 & --&2.76 &--&--&0.0102\\
& 4 & 0.0160& 0.0340&4&26021 & --& 0.58&--&--&0.0160\\
& 5 & 0.0135 & 0.0314 &5&125555 & --& 2.61&--&--&0.0135\\
& 10 & 0.0119  & 0.0290 &10&35025 & --& 1.08&--&--& 0.0119\\
& 15 & 0.0117 & 0.0286 &15& 4705& --& 0.42&--&--&0.0117\\
& 20 & 0.0104  &0.0282 &20&1102 & --&0.44 &--&--&0.0104\\
& 25 & 0.0104 & 0.0282 & 24& 909& --&0.48 &--&--&0.0104\\
& 30 & 0.0104 & 0.0282 &28 &545 & --&0.50 &--&--&0.0104\\
& 35 & 0.0107 & 0.0282& 32 & 14 & --& 0.02 &--&--&0.0107\\
& 40 & 0.0104 & 0.0282 & 33 & 13 & --&0.02 &--&--&0.0104\\
& 45 & 0.0114 &0.0283 &30 &13 & --&0.03 &--&--&0.0114\\
& 50 & 0.0104 & 0.0282 & 33& 13& --& 0.03&--&--&0.0104\\
& 55 & 0.0104  & 0.0282 & 33 & 13& --&0.03 &--&--&0.0104\\
& 60 & 0.0133 &0.0289 & 27 &13 & --&0.03 &--&--&0.0133\\
& 65 & 0.0145 & 0.0294 & 28 &14 & --&0.05 &--&--&0.0145\\
& 70 & 0.0105  & 0.0282 & 33 & 13 & --& 0.02 &--&--&0.0105\\
& 75 & 0.0105 & 0.0282 & 33 & 13 & --& 0.02 &--&--&0.0105\\
& 80 &  0.0105 & 0.0282 & 33 & 13 & --&  0.02 &--&--&0.0105\\
& 85 & 0.0104  & 0.0282 & 33 & 13 & --& 0.02&--&--&0.0104\\
& 89 & 0.0101 & 0.0282& 34 &13 & --&0.09 &--&--&0.0101\\
\hline
\end{tabular}
\end{center}
\end{table}

\begin{table}[H]
	\caption{Performance of Algorithm PSPGD and CPLEX for problem Port4.} \label{Table:port4_b}
	\begin{center}
		\begin{tabular}{|ccccccccccc|}\hline
			Problem&$\alpha$ & $v^{\top}x^*$ &$\sqrt{(x^*)^{\top}Qx^*}$ & $\|x^*\|_0$ & Iter & Iter-SPG & Time & $\tau$  & fcnt & $\rho$\\ \hline \hline
				PSPGD&  1 & 0.0343 &  0.0983 &         1    &          2  &  11 & 1.0960  &  0.0891 &    32 & 0.0095 \\
	 &2 & 0.0368 &  0.1084 &         1    &          2  &   5 & 0.2783  &  0.8903 &   18 & 0.0091 \\
	&3 & 0.0368 &  0.1084 &         1    &          2  &   5 & 0.2805  &  0.8899  &   18 & 0.0075 \\
	&4 & 0.0368&  0.1084 &         1    &          2  &  26 & 0.4912  &  0.8908  &   74 & 0.0108 \\
	&5 & 0.0368 &  0.1084 &         1    &          2  &   5 & 0.2296  &  0.8905  &   18 & 0.0101 \\
	&10 & 0.0368 &  0.1084 &         1    &          2  &   4 & 0.2047  &  0.8890  &   11 & 0.0050 \\
&15 & 0.0368 &  0.1084 &         1    &          2  &   5 & 0.2230  &  1.7778 &   17 & 0.0047 \\
  &20 & 0.0368 &  0.1084 &         1    &          2  &   4 & 0.3299  &  3.7049  &   12 & 0.0048 \\
 & 25 & 0.0206 &  0.0365 &        10    &          3  &  72 & 0.9029  &  0.4027 &  118 & 0.0053 \\
& 30 & 0.0194 &  0.0346 &        14    &          3  & 107 & 2.0201  &  0.5036 &  194 & 0.0060 \\
&  35 & 0.0197 &  0.0349 &        13    &          3  &  74 & 0.9154  &  0.5033  &  133 & 0.0060 \\
&40 & 0.0178 &  0.0319 &        20    &          3  &  96 & 1.7136  &  0.8921  &  144 & 0.0055 \\
& 45 & 0.0127 &  0.0371 &         4    &          3  &  51 & 0.9430  &  0.0202&   54 & 0.0091 \\
& 50 & 0.0175 &  0.0410 &         3    &          3  &  51 & 0.9916  &  0.0202 &   67 & 0.0091 \\	
& 55 & 0.0109 &  0.0292 &         6    &          2  &  24 & 0.4256  &  0.0101  &   34 & 0.0106 \\
&60 & 0.0194 &  0.0346 &        14    &          6  & 213 & 3.8296  &  0.9062  &  417 & 0.0063 \\
&     65 & 0.0191 &  0.0339 &        16    &         20  & 979 & 17.236  &  1.1020  & 2148 & 0.0061 \\
&70 & 0.0132 &  0.0344 &         5    &          2  &  20 & 0.3864  &  0.0301 &   28 & 0.0010 \\
& 75 & 0.0252 &  0.0481 &         7    &         13  & 561 & 14.115  &  0.9461 & 1052 & 0.0067 \\
			&80 & 0.0138 &  0.0365 &         5    &          2  &  27 & 2.1716  &  0.0888 &   35 & 0.0080 \\
			&85 & 0.0142 &  0.0364 &         6    &          2  &  17 & 1.5391  &  0.0888  &   25 & 0.0075 \\
			&90 & 0.0086 &  0.0250 &        14    &          2  &  23 & 1.1892  &  0.0888 &   25 & 0.0073 \\
			&95 & 0.0089 &  0.0231 &        18    &          2  &  28 & 1.1098  &  0.0888  &   30 & 0.0080 \\
			&98 & 0.0098 &  0.0223 &        38    &          2  &  20 & 0.9998  &  0.0888  &   22 & 0.0098 \\				
			\hline   \hline
			CPLEX
			& 1  & 0.0115 & 0.0462&1 &806 & --&0.09 &--&--&0.0095\\
			& 2  & 0.0095 & 0.0350&2 &29245 & --&0.48 &--&--&0.0091\\
			& 3  & 0.0081 & 0.0300&3 &506050 & --& 3.63&--&--&0.0075\\
			& 4  & 0.0108 & 0.0287& 4& 1750081& --&14.89 &--&--&0.0108\\
			& 5  & 0.0101& 0.0266 &5 & 2497651 & --&20.86&--&--&0.0101\\
			& 10  & 0.0070 & 0.0231 & 10& 698137& --&6.38&--&--&0.0050\\
	& 15  & 0.0077 & 0.0223 & 15&8163 & --&0.33&--&--&0.0047\\
			& 20  & 0.0075 &0.0221 &20 &11669&--&0.33&--&--&0.0048\\
			& 25 &0.0074  &0.0221 &25 &1062& --&0.19&--&--&0.0053\\
			& 30  &0.0078&0.0221 &28 &540 & --&0.17 &--&--&0.0060\\
			& 35  &0.0076&0.0221 &35 &74 & --&0.09 &--&--&0.0060\\
			& 40  & 0.0077 &0.0220 &38 &14 & --&0.02 &--&--&0.0055\\
			& 45 & 0.0091 &0.0222 & 39 &14 & --&0.02 &--&--&0.0091\\
			& 50 & 0.0091 & 0.0222& 39&14 & --&0.00 &--&--&0.0091\\
			& 55 & 0.0106&0.0226&35&13 & --&0.01&--&--&0.0106\\
			& 60 & 0.0077 &0.0220 &38 &14 & --&0.03 &--&--&0.0063\\
			& 65 & 0.0077 &0.0220 &38 &14 & --&0.02 &--&--&0.0061\\
			& 70 & 0.0077 & 0.0220 &38 &14 & --& 0.03&--&--&0.0010\\
			& 75 & 0.0077 &0.0220 &38 &14 & --&0.06 &--&--&0.0067\\
			& 80 & 0.0080  &0.0220 &38 &14 & --&0.02 &--&--&0.0080\\
			& 85 & 0.0077 &0.0220 &38 &14 & --&0.01 &--&--&0.0075\\
			& 90 & 0.0077 &0.0220 &38 &14 & --&0.02 &--&--&0.0073\\
		    & 95 & 0.0080 & 0.0220&38 &14 & --&0.02 &--&--&0.0080\\
		    & 98 & 0.0098 & 0.0223 & 38 & 14& --& 0.03&--&--&0.0098\\
			\hline
		\end{tabular}
	\end{center}
\end{table}

	\begin{table}[H]
		\caption{Performance of Algorithm PSPGD and CPLEX for problem Port5.} \label{Table:port5_b}
		\begin{center}
			\begin{tabular}{|ccccccccccc|}\hline
				Problem&$\alpha$ & $v^{\top}x^*$ &$\sqrt{(x^*)^{\top}Qx^*}$ & $\|x^*\|_0$ & Iter & Iter-SPG & Time & $\tau$  & fcnt & $\rho$\\ \hline \hline
				PSPGD	& 2 & 0.0161 &  0.1081 &         2    &          2  & 102 & 7.4534  &  0.1010 &  356 & 7.7209e-06 \\

      & 3 & 0.0037 &  0.0538 &         3    &         20  & 1020 & 77.895  &  0.0100 & 1342 & 7.8605e-06 \\			
&			  4 & 0.0034 &  0.0500 &         3    &         20  & 1020 & 80.528  &  0.0100  & 1333 & 7.8656e-06 \\
		    &5 & 0.0009 &  0.0388 &         4    &         12  & 612 & 43.506  &  0.0100 &  772 & 1.2051e-05 \\		
		    &  10 & 0.0006 &  0.0355 &         7    &          4  & 204 & 16.183  &  0.0400 &    482 & 1.1788e-05 \\
			     & 15 & 0.0036 &  0.0377 &         8    &          5  & 255 & 19.572  &  0.2000 &  658 & 7.2419e-06 \\						
     & 20 & 0.0059&  0.0399 &         9    &         11  & 561 & 39.421  &  0.9000  & 1633 & 5.9825e-06 \\
	      &25 & 0.0058 &  0.0399 &        10    &          9  & 459 & 19.429  &  1.0000 & 1461 & 6.4813e-06 \\						
			&			  30 & 0.0003 &  0.0349 &        10    &          8  & 325 & 21.612  &  0.0001 &  463 & 7.8618e-06 \\
			     &225 & 0.0003 &  0.0349 &       12     &          2  &  10 & 1.4476 &  0.9051 &   12 & 1.1926e-05 \\
					
			\hline   \hline
CPLEX 	&2 &0.0058 &0.0439 &2&1964&  -- &  0.48 &  -- &   -- & 7.7209e-06 \\
&3 & 0.0027&0.0391&3&1034&  --&0.41  &  --  &   -- & 7.8605e-06 \\
&4 & 0.0009&0.0367 &4&386&  --&0.61 &  -- &  -- & 7.8656e-06 \\
&	5 & 0.0003 &0.0356&5&132&-- &0.33 &  --  &   -- & 1.2051e-05 \\
&10 &0.0003 &0.0349&10&19&   -- &0.27 &  --  &   --& 1.1788e-05 \\
&15 &0.0003 &0.0349&12&17&   -- &0.22   &  --&   -- & 7.2419e-06 \\
&	 20 & 0.0003&0.0349&12&17&   -- &0.11   & -- &   -- & 5.9825e-06 \\
&	25 & 0.0003&0.0349&12&17&   -- &0.26  &  --  &   -- & 6.4813e-06
\\
&	30 & 0.0003&0.0349&12&17& -- & 0.09  &  -- &  -- & 7.8618e-06 \\
		    & 225 &0.0003 &0.0349& 12 &17 & --&0.11 &--&--&1.1926e-05\\
			\hline
		\end{tabular}
	\end{center}
\end{table}

\begin{table}[H]
				\caption{Performance of Algorithm PSPGD for problem Port6.} \label{Table:port6}
				\begin{center}
					\begin{tabular}{|ccccccccccc|}\hline
						Problem&$\alpha$ & $v^{\top}x^*$ &$\sqrt{(x^*)^{\top}Qx^*}$ & $\|x^*\|_0$ & Iter & Iter-SPG & Time & $\tau$  & fcnt & $\rho$\\ \hline \hline
	
					PSPGD &
					
1 & 0.0364 &  0.0909 &         1    &          2  &  15 & 1.3413  &  1.7475 &   53 & 1.3020e-03 \\
&4 & 0.0969 &  0.1578 &         3    &          2  &  24 & 1.3470  &  1.7480 &  108 & 1.3020e-03 \\
&5 & 0.0969 &  0.1578 &         3    &          2  &  31 & 1.5071  &  1.7480 &   143 & 1.3020e-03 \\
  &10 & 0.0709 &  0.1845 &         2    &          2  &  11 & 0.8574  &  52.424 &    50 & 1.3020e-03 \\
&15 & 0.0364 &  0.0909 &         1    &          2  &   4 & 0.6542  &  52.425 &    16 & 1.3020e-03 \\
&20 & 0.0364 &  0.0909 &         1    &          2  &   4 & 0.6266  &  52.425 &    16 & 1.3020e-03 \\
&25 & 0.0364 &  0.0909 &         1    &          2  &   4 & 0.5857  &  52.424 &    16 & 1.3020e-03 \\
&	30 & 0.0364 &  0.0909 &         1    &          2  &   4 & 0.5743  &  52.424 &   16 & 1.3020e-03 \\
			
&35 & 0.0548 &  0.2302 &    3    &          4  & 113 & 2.6844  &  34.984 &   434 & 1.3020e-03 \\
					
&38 & 0.0364 &  0.0909 &         1    &          2  &   4 & 1.1503  &  52.424 &   16 & 1.3020e-03 \\
					
&235 & 0.0709 &  0.1845 &         2    &          2  &   8 & 0.9253  &  1.7471 &     13 & 1.3020e-03 \\
&245 & 0.0709 &  0.1845 &         2    &          2  &   8 & 0.9639  &  1.7471 &    13 & 1.3020e-03 \\
&255 & 0.0204 &  0.3015 &         2    &          2  &  29 & 1.3874  &  1.7467 &   106 & 1.3020e-03 \\
&265 & 0.0204 &  0.3015 &         2    &          2  &  29 & 1.4178  &  1.7467 &    106 & 1.3020e-03 \\
&275 & 0.0204 &  0.3015 &         2    &          2  &  29 & 1.3963  &  1.7467 &   106 & 1.3020e-03 \\
&285 & 0.0204 &  0.3015 &         2    &          2  &  29 & 1.5191  &  1.7467 &    106 & 1.3020e-03 \\
&295 & 0.0204 &  0.3015 &         2    &          2  &  11 & 0.9662  &  1.7467 &     16 & 1.3020e-03 \\
&305 & 0.0204 &  0.3015 &         2    &          2  &  13 & 1.1448  &  1.7467 &     27 & 1.3020e-03 \\
&315 & 0.0204 &  0.3015 &         2    &          2  &  12 & 1.0017  &  1.7467 &     17 & 1.3020e-03 \\
&325 & 0.0204 &  0.3015 &         2    &          2  &  12 & 1.0606  &  1.7467 &     17 & 1.3020e-03 \\

&525 & 0.0139 &  0.0278 &        19    &          4  & 102 & 10.025  &  1.7475 &   153 & 1.3060e-03 \\
&535 & 0.0137 &  0.0276&        19    &          7  & 255 & 22.419  &  1.7492 &    419 & 1.3059e-03 \\
&545 & 0.0139 &  0.0278 &        19    &          5  & 153 & 14.205  &  1.7477 &    207 & 1.3059e-03 \\
&555 & 0.0146 &  0.0285 &        20    &          5  & 153 & 14.119  &  1.7475 &   232 & 1.3059e-03 \\
&565 & 0.0183 &  0.0340 &        18    &          4  & 102 & 10.159  &  1.7470 &    116 & 1.3058e-03 \\
&575 & 0.0125 &  0.0263 &        20    &          3  &  51 & 5.9195  &  1.7467 &     54 & 1.3058e-03 \\
&585 & 0.0084 &  0.0244 &        23    &          3  &  61 & 6.6879  &  1.7467 &    81 & 1.3057e-03 \\
&595 & 0.0046 &  0.0226 &        24    &          2  &  46 & 5.1943  &  1.7466 &     48 & 1.3057e-03 		\\
&600 & 0.0015 &  0.0207 &      38     &          2  &  40 & 2.0572  &  1.7466 &   43 & 1.3057e-03 \\

					\hline
					
							\end{tabular}
				\end{center}
			\end{table}

\end{document}